\documentclass[a4paper,12pt,reqno]{amsart}
\usepackage{amsfonts}
\usepackage{amsmath}
\usepackage{amssymb}
\usepackage[a4paper]{geometry}
\usepackage{mathrsfs}
\usepackage{csquotes}
\usepackage{enumitem}
\usepackage{dirtytalk}
\usepackage{comment}
\usepackage{nccmath}

\usepackage[colorlinks]{hyperref}
\renewcommand\eqref[1]{(\ref{#1})} %Need with hyperref
\numberwithin{equation}{section}
\theoremstyle{plain}
\newtheorem{thm}{Theorem}[section]
\newtheorem{prop}[thm]{Proposition}
\newtheorem{cor}[thm]{Corollary}
\newtheorem{lem}[thm]{Lemma}

\theoremstyle{definition}
\newtheorem{defn}[thm]{Definition}
\newtheorem{rem}[thm]{Remark}

\newcommand\norm[1]{\left\lVert#1\right\rVert}

\def\e[#1]{{\textrm{e}}^{#1}}

\begin{document}

   \title[General weighted Hardy-Poincar\'e type inequalities]
 {Sharp remainder terms and stability of weighted Hardy-Poincar\'e and Heisenberg-Pauli-Weyl inequalities related to the Baouendi-Grushin operator} 

   \author[Y. Shaimerdenov]{Yerkin Shaimerdenov}
\address{
  Yerkin Shaimerdenov:
 \endgraf
   SDU University, Kaskelen, Kazakhstan
  \endgraf
  and 
  \endgraf
  Department of Mathematics: Analysis, Logic and Discrete Mathematics 
  \endgraf
  Ghent University, Belgium
  \endgraf
  and
  \endgraf
  Institute of Mathematics and Mathematical Modeling, Kazakhstan
   \endgraf
  {\it E-mail address} {\rm
yerkin.shaimerdenov@sdu.edu.kz}
  }

\author[N. Yessirkegenov]{Nurgissa Yessirkegenov}
\address{
  Nurgissa Yessirkegenov:
  \endgraf
  KIMEP University, Almaty, Kazakhstan
     \endgraf
  {\it E-mail address} {\rm nurgissa.yessirkegenov@gmail.com}
  }

  \author[A. Zhangirbayev]{Amir Zhangirbayev}
\address{
  Amir Zhangirbayev:
 \endgraf
   SDU University, Kaskelen, Kazakhstan
  \endgraf
  and 
  \endgraf
  Institute of Mathematics and Mathematical Modeling, Kazakhstan
   \endgraf
  {\it E-mail address} {\rm
amir.zhangirbayev@gmail.com}
  }

\thanks{This research is funded by the Committee of Science of the Ministry of Science and Higher Education of the Republic of Kazakhstan (Grant No. AP23490970).}


     \begin{abstract} In this paper, we obtain sharp remainder terms for the Hardy-Poincar\'e inequalities with general non-radial weights in the setting of Baouendi-Grushin vector fields (see Theorem \ref{main thm}). It is worth emphasizing that all of our results are new both in the Baouendi-Grushin and standard Euclidean settings. The method employed allows us to not only unify, but also improve the results of Kombe and Yener \cite{KY18} for any $1<p<\infty$ while holding true for complex-valued functions and providing explicit constants (Corollary \ref{cor rec KY}). As a result, we are able to obtain sharp remainder terms to many known weighted Hardy-type inequalities (see Section \ref{sec sharp remainders}). Aside from weighted Hardy-type inequalities, we also recover a sharp remainder formula for the $L^{p}$-Poincar\'e inequality (Corollary \ref{cor poin}). In the special case of radial weights, we are naturally able to introduce the notion of Baouendi-Grushin $p$-Bessel pairs (see Definition \ref{p-bessel pair}). Furthermore, we apply the technique to establish the sharp remainder term of the Heisenberg-Pauli-Weyl inequality in $L^{p}$ for $1<p<\infty$ (Corollary \ref{cor hpw sharp}), which includes the sharp constant. This makes it possible to obtain the $L^{p}$-analogue for $2\leq p < n$ (Theorem \ref{stability hpw}) of a stability result by Cazacu, Flynn, Lam and Lu \cite{cazacu2024caffarelli}. Lastly, as another application, the non-existence of positive solutions to nonlinear parabolic partial differential equations is investigated (Theorem \ref{parabolic eq}).

     \end{abstract}
     \maketitle     

\section{Introduction}

\allowdisplaybreaks

The classical Hardy inequality is a cornerstone in the analysis of partial differential equations and functional inequalities. Its multidimensional (in $L^{p}$) form is stated as follows: let $\Omega\subseteq\mathbb{R}^{n}$, $1 < p < n$, then, for all $f \in C_0^{\infty}(\Omega)$, the following inequality holds:
\begin{equation} \label{class_multdim_hardy}
\int_{\Omega} |\nabla f|^p dz \geq \left( \frac{n - p}{p} \right)^p \int_{\Omega} \frac{|f|^p}{|z|^p} dz,
\end{equation}
where the constant $\left( \frac{n - p}{p} \right)^p$ is sharp but not attained.

One of the well-known applications of the Hardy inequality of the form (\ref{class_multdim_hardy}) lies in the study of solutions to partial differential equations. For example, in \cite{baras1984heat}, Baras and Goldstein found that the existence of solutions to the heat equation with singular potentials is connected to the extremizers and sharp constants of the Hardy inequality. Later, Goldstein and Zhang generalized these results to heat equations on the Heisenberg group \cite{goldstein2001degenerate}, where they used the Hardy inequality on the Heisenberg group by Garofalo and Lanconelli \cite{garofalo1990frequency}. In the meantime, Brezis and V\'azques \cite{brezis1997blow} studied blow-up solutions of some nonlinear elliptic problems. There, they obtained a characterization of extremal solutions to a certain semilinear elliptic equation. However, in order to apply this characterization to standard examples, they needed an improved Hardy inequality, which they successfully derived \cite[Theorem 4.1]{brezis1997blow}. Brezis and V\'azques also posted an open problem \cite[Problem 2 of Section 8]{brezis1997blow}, which states whether the improved Hardy inequality \cite[Theorem 4.1]{brezis1997blow} can be further improved in a certain context. 

The findings of Baras, Goldstein, Zhang, Brezis and V\'azques received enormous attention and led to the development of Hardy inequalities, their improvements and generalizations (e.g. \cite{brezis1997hardy, vazquez2000hardy, chaudhuri2002improved, filippas2002optimizing, gazzola2004hardy, barbatis2004unified, ghoussoub2008best, cazacu2013improved}). We also refer to standard monographs \cite{KMP06, KMP07, maz2011sobolev,  ghoussoub2013functional, BEL15, persson2017weighted, ruzhansky2019hardy} along with open problems posed by Maz'ya on Hardy inequalities and many other topics \cite{maz2018seventy}.

A specific line of development in Hardy inequalities that we aim to focus on originates with the results by Ghoussoub and Moradifam \cite{ghoussoub2011bessel}, where they gave a sufficient and necessary condition on a pair of positive radial functions $v(|z|)$ and $w(|z|)$ so that, we have
\begin{align*}
\int_{B_{R}}v(|z|)|\nabla f|^{2}dz\geq\int_{B_{R}}w(|z|)|f|^{2}dz
\end{align*}
for all $f\in C^{\infty}_{0}(B_{R})$. The results of Ghoussoub and Moradifam were later extended and further developed (see e.g. \cite{moradifam2011optimal, skrzypczak2014hardy, lam2018hardy, ruzhansky2018extended, lam2019factorizations, duy2020improved, lam2020geometric, flynn2021sharp, duy2022p, ruzhansky2024hardy, GJR24, flynn2025p}). See also \cite{ghoussoub2013functional} for many examples and properties of Bessel pairs. 

After the contributions of Ghoussoub, Moradifam and others, considerable effort has been devoted to the study of non-radial weights, leading to xsignificant advances in this direction. For example, in \cite{frank2008non}, Frank and Seiringer developed non-linear and non-local versions of the ground state representation that allow to refine, generalize
and unify Hardy inequalities without the requirement on radial weights: assume that $\phi$ is a positive weak solution of the weighted $p$-Laplace equation 
\begin{align}\label{frank}
-\nabla \cdot\left(v\left|\nabla \phi\right|^{p-2} \nabla\phi\right) = w\phi^{p-1}.
\end{align}
Then, one has
\begin{align*}
\int_{\mathbb{R}^n}v|\nabla f|^{p}dz\geq\int_{\mathbb{R}^n}w|f|^{p}dz.
\end{align*}
By utilizing the convexity inequality for $p\geq2$, one can obtain the same inequality with a positive remainder term \cite[Formula (2.12)]{frank2008non}. In the same spirit, Goldstein, Kombe and Yener in \cite{goldstein2017unified}, applied this approach to establish general weighted Hardy inequalities on a Carnot group. In contrast to the equality in (\ref{frank}), they had an inequality, thereby relaxing the condition and giving them a broader class of Hardy inequalities. This idea subsequently found applications in many different contexts (e.g. \cite{goldstein2018general, KY18, yener2018general, yener2018several, ruzhansky2019weighted, cazacu2021method}). We also refer to \cite{hoffmann2015hardy} to Hardy inequalities with homogeneous weights.

However, we are mostly interested in \cite{KY18}. Before discussing the results of the paper, we first would like to introduce the Baouendi-Grushin setting and its relation to Hardy inequalities, which will be central to our analysis. 

We consider the sub-Riemannian space $\mathbb{R}^{m+k}=\{z=(x,y):x\in\mathbb{R}^{m}, y\in\mathbb{R}^{k}\}$ defined by the Baouendi-Grushin vector fields 
\begin{align}\label{grushin vector fields}
X_i = \frac{\partial}{\partial x_i}, \quad i = 1, \ldots, m, \qquad
Y_j = |x|^{\gamma} \frac{\partial}{\partial y_j}, \quad j = 1,\ldots,k.
\end{align}
Here $\gamma\geq0$ is any nonnegative real number with $|x|$ representing the usual Euclidean  norm of $x\in\mathbb{R}^{m}$. After being initially introduced by Baouendi \cite{baouendi1967classe} and later by Grushin \cite{gruvsin1970class, gruvsin1971class}, Garofalo was the first to obtain the Hardy inequality in this setting \cite[Corollary 4.1]{Gar93}:
\begin{align}\label{Gar}
\int_{\mathbb{R}^n}\left(|\nabla_{x}f|^{2}+|x|^{2\gamma}|\nabla_{y}f|^{2}\right) dz\geq \left(\frac{Q-2}{2}\right)^2 \int_{\mathbb{R}^n}\frac{|x|^{2 \gamma}}{|x|^{2\gamma+2}+(1+\gamma)^{2}|y|^{2}}|f|^2 dz,
\end{align}
where $n=m+k$, $Q=m+(1+\gamma)k$ and $f \in C_0^{\infty}\left(\mathbb{R}^m \times \mathbb{R}^k \backslash\{(0,0)\}\right)$, $\nabla_{x}f$ is the gradient of $f$ in the variable $x$ and $\nabla_{y}f$ is in the variable $y$. In the special case of $\gamma=0$, the inequality (\ref{Gar}) reduces to the classical multidimensional $L^{2}$-Hardy inequality (i.e. $p=2$ of (\ref{class_multdim_hardy})).

Let us now give some additional facts on the Baouendi-Grushin framework. The corresponding gradient given by the vector fields (\ref{grushin vector fields}) is defined by
\begin{align*}
\nabla_{\gamma}=\left(X_{1},\ldots,X_{m},Y_{1},\ldots,Y_{k}\right)=(\nabla_{x},|x|^{\gamma}\nabla_{y}).
\end{align*}
The Baouendi-Grushin operator is defined as a differential operator on $\mathbb{R}^{m+k}$ by
\begin{align*}
\Delta_{\gamma}=\sum_{i=1}^{m}X^{2}_{i}+\sum_{j=1}^{k}Y^{2}_{j}=\Delta_{x}+|x|^{2\gamma}\Delta_{y}=\nabla_{\gamma}\cdot\nabla_{\gamma}.
\end{align*}
The operator in itself has many interesting properties. For example, it is not uniformly elliptic in the space $\mathbb{R}^{m+k}$ (due to the degeneracy on the subspace $\{0\}\times\mathbb{R}^{k}$). It is also sub-elliptic and for non-negative even integers $\gamma$ satisfies the H\"ormander condition, meaning it can be expressed as a sum of squares of vector fields that form a Lie algebra of full rank at any point in $\mathbb{R}^{m+k}$. We also have the anisotropic dilation defined by
\begin{align*} 
\delta_{a}(x,y)=(ax,a^{1+\gamma}y),
\end{align*}
where $a>0$. The homogeneous dimension with respect to the dilation is $Q=m+(1+\gamma)k$. The formula of the change of variables for the Lebesgue measure gives
\begin{align*}
d \circ \delta_a(x, y)=a^Q d x d y.
\end{align*}
The corresponding distance function $\rho(z)$ is given by
\begin{align*}
\rho=\rho(z):=\left(|x|^{2(1+\gamma)}+(1+\gamma)^2|y|^2\right)^{\frac{1}{2(1+\gamma)}}.
\end{align*}
Direct computations give us the following formulas
\begin{align*}
\nabla_{\gamma}\rho=\left(\frac{|x|^{2\gamma}x}{\rho^{2\gamma+1}},\frac{(1+\gamma)|x|^{\gamma}y}{\rho^{2\gamma+1}}\right) \quad \text{and} \quad |\nabla_{\gamma}\rho|=\frac{|x|^\gamma}{\rho^\gamma}.
\end{align*}
Also, we have
\begin{align*}
\nabla_{\gamma}\cdot\left(\rho^{c}|x|^{s}\nabla_{\gamma}\rho\right)=(Q+c+s-1)\frac{|x|^{2\gamma+s}}{\rho^{2\gamma+1-c}}
\end{align*}
for $c,s\in\mathbb{R}$. Additionally, we define the so-called $p$-Grushin operator similarly to how the classical Baouendi-Grushin operator is derived from the Laplacian:
\begin{align*}
\Delta_{\gamma,p}\phi=\nabla_{\gamma}\cdot\left(|\nabla_{\gamma}\phi|^{p-2}\nabla_{\gamma}\phi\right)
\end{align*}
for some complex-valued $\phi$ on $\Omega\subseteq\mathbb{R}^{m+k}$, which gives us the following formula
\begin{align*}
\Delta_{\gamma,p} \phi(\rho)
= \frac{|x|^{\gamma p}}{\rho^{\gamma p}} \, \left| \phi'(\rho) \right|^{p-2}
\left[ (p-1)\phi''(\rho) + (Q-1)\frac{\phi'(\rho)}{\rho} \right], 
\quad 1<p<\infty
\end{align*}
for all $z\in\mathbb{R}^{m+k} \backslash\{(0,0)\}$. Lastly, we define $B_{R}=\{z\in\mathbb{R}^{m+k}:\rho(z)<R\}$ as the $\rho$-ball centered at the origin with radius $R>0$.

Returning to Hardy inequalities, we continue from the result of Garofalo \cite[Corollary 4.1]{Gar93}. After that, there has been an intense development. In \cite{D'A04}, D'Ambrosio obtained weighted $L^{p}$-versions of (\ref{Gar}) as well as the $L^{p}$-Poincar\'e inequality contained in a slab with a particular assumption on the domain. Then, Niu, Chen and Han \cite{NCH04} employed  the Picone identity from \cite{zhang2003picone} to obtain many different Hardy inequalities on $B_{R}$, $\mathbb{R}^{m+k}\backslash B_{R}$ and on $\mathbb{R}^{m+k}$. Later, D'Ambrosio \cite{D'A05} further developed and unified the results from \cite{D'A04} by establishing weighted $L^{p}$-Hardy inequalities, this time for general quasilinear second-order differential operators (which include the Baouendi-Grushin operator). Over the next years, Hardy inequalities for Baouendi-Grushin vector fields have remained a huge topic of interest and we refer the reader to many interesting works in this direction \cite{SJ12, Kom15, yang2015improved, LRY19, suragan2023sharp, d2024weighted, d2024unified, GJR24, yessirkegenov2025refined, banerjee2025extremizer}.

Now the result that we are mostly focused on is from \cite{KY18}, where Kombe and Yener obtained the following general weighted $L^{p}$-Hardy type inequality: let $v \in C^1\left(\mathbb{R}^n\right)$ and $w \in L_{\text {loc }}^1\left(\mathbb{R}^n\right)$ be nonnegative functions and $\phi \in$ $C^{\infty}\left(\mathbb{R}^n\right)$ be a positive function satisfying the differential inequality
\begin{align}\label{condt kombe}
-\nabla_\gamma \cdot\left(v\left|\nabla_\gamma \phi\right|^{p-2} \nabla_\gamma\phi\right) \geq w\phi^{p-1}
\end{align}
a.e. in $\Omega$. There is a positive number $c_p=c(p)$ such that; for $p \geq 2$, then
\begin{align}\label{kom1}
\int_{\mathbb{R}^n} v\left|\nabla_\gamma f\right|^p dz \geq \int_{\mathbb{R}^n} w|f|^p dz+c_p \int_{\mathbb{R}^n} v\left|\nabla_\gamma \left(\frac{f}{\phi}\right)\right|^p \phi^p dz
\end{align}
and for $1<p<2$, then
\begin{align}\label{kom2}
\int_{\mathbb{R}^n} v\left|\nabla_\gamma f\right|^p dz \geq \int_{\mathbb{R}^n} w|f|^p dz+c_p \int_{\mathbb{R}^n} \frac{v\left|\nabla_\gamma \left(\frac{f}{\phi}\right)\right|^2 \phi^2}{\left(\left|\frac{f}{\phi} \nabla_\gamma \phi\right|+\left|\nabla_\gamma \left(\frac{f}{\phi}\right)\right| \phi\right)^{2-p}} dz
\end{align}
for all real-valued functions $f \in C_0^{\infty}\left(\mathbb{R}^n\right)$.

The results, (\ref{kom1}) and (\ref{kom2}), in themselves are very interesting as they recover and improve most known Hardy inequalities (e.g. \cite[Theorem 2.1]{NCH04}, \cite[Formula (3.1) of Theorem 3.1]{D'A04}, \cite[Formula (3.11) of Theorem 3.5]{D'A05} and etc.) with radial and non-radial weights for all real-valued functions $f \in C_0^{\infty}\left(\mathbb{R}^n\right)$. Not only that, they also derive new inequalities in the context of the Baouendi-Grushin operator (e.g. \cite[Theorem 5.3]{skrzypczak2013hardy}, \cite[Formula (7)]{ghoussoub2011bessel}).

The main goal of this paper is to provide a way to obtain sharp remainder terms (identities) and explicit constants in (\ref{kom1}) and (\ref{kom2}) without the real-valued assumption. Also, instead of having two different inequalities for each case of $p$, it would be good to have one unifying inequality that contains them both. Finally, the condition (\ref{condt kombe}) seems to work when $v=1$ and $\phi:=\text{the first eigenfunction of $-\Delta_{\gamma, p}$}$, which should give the $L^{p}$-Poincare inequality. However, the first eigenfunction of $-\Delta_{\gamma, p}$ is not smooth in general. So, it is interesting to see whether the condition on $\phi$ can be made less restrictive.

In this paper, we do all of the above. More precisely, we have the following result: let $1<p<\infty$ and let $\Omega\subseteq\mathbb{R}^{m+k}$ be an open set such that the integrals below make sense. Let $v,w\in L^{1}_{loc}(\Omega)$ be nonnegative functions with $\phi>0$ satisfying the following in the distributional sense in $\Omega$:
\begin{align}\label{intro main condition}
-\nabla_{\gamma} \cdot \left( v \left|\nabla_{\gamma} \phi\right|^{p-2} \nabla_{\gamma} \phi \right)\geq w\phi^{p-1}.
\end{align}
Then, for all complex-valued $f\in C^{\infty}_{0}(\Omega)$, we have 
\begin{align}\label{intro main result}
\int_{\Omega}^{}v|\nabla_{\gamma}f|^{p}dz \geq \int_{\Omega}^{}w|f|^{p}dz + \int_{\Omega}^{}vC_{p}\left( \nabla_{\gamma}f,\phi\nabla_{\gamma}\left( \frac{f}{\phi} \right)  \right) dz.
\end{align}
Equality, in (\ref{intro main result}), holds when $\phi$ is a positive distributional solution to (\ref{intro main condition}) with an equality. 

Following the approach of~\cite{KY18}, one can derive arbitrarily many sharp remainder terms for weighted \( L^{p} \)-Hardy type inequalities, provided there exists positive \( w \) for some given \( v \) and \( \phi \) such that ~\eqref{intro main condition} holds with an equality. To illustrate, we can choose 
\begin{align*}
v=|x|^{\beta-\gamma p}\rho^{(1+\gamma)p-\alpha} \quad \text{and} \quad \phi=\rho^{-\frac{Q+\beta -\alpha}{p}},
\end{align*}
which from (\ref{intro main condition}) gives $w=\left(\frac{Q+\beta-\alpha}{p}\right)^{p}\frac{|x|^{\beta}}{\rho^{\alpha}}$ and, hence, the sharp remainder term of \cite[Formula (3.1) of Theorem 3.1]{D'A04}: let $1<p<\infty$, $m,k\geq1$ and $\alpha,\beta \in\mathbb{R}$ such that $Q\geq\alpha-\beta$. Then, for all complex-valued $f\in C^{\infty}_{0}(\mathbb{R}^{m+k}\backslash\{(0,0)\})$, we have
\begin{multline}\label{dam org eq intro}
\int_{\mathbb{R}^n}^{}|x|^{\beta-\gamma p}\rho^{(1+\gamma)p-\alpha}|\nabla_{\gamma}f|^{p}dz=\left( \frac{Q+\beta-\alpha}{p} \right) ^{p}\int_{\mathbb{R}^n}^{}\frac{|x|^{\beta}}{\rho^{\alpha}}|f|^{p}dz\\+\int_{\mathbb{R}^n}^{}|x|^{\beta-\gamma p}\rho^{(1+\gamma)p-\alpha}C_{p}\left( \nabla_{\gamma}f,\rho^{-\frac{Q+\beta-\alpha}{p}}\nabla_{\gamma}\left(  \frac{f}{\rho^{-\frac{Q+\beta-\alpha}{p}}}\right)  \right)dz.
\end{multline}
Identity (\ref{dam org eq intro}), under an appropriate choice of parameters, immediately gives the sharp remainder term of \cite[Formula (3.2) of Theorem 3.1]{D'A04} and (\ref{Gar}). Additionally, (\ref{dam org eq intro}), in the Euclidean case, recovers \cite[Lemma 3.3]{cazacu2024hardy}.

Here, the result (\ref{intro main result}) is cleanly written as one inequality that (due to \cite[Step 3 of Proof of Lemma 3.4]{cazacu2024hardy} and \cite[Lemma 2.2]{CT24}) contains both (\ref{kom1}) and (\ref{kom2}) at once with explicit constants.

Moreover, we are able to relax the condition on $\phi$ in (\ref{condt kombe}) by viewing the condition (\ref{condt kombe}) as a $p$-Grushin equation with a positive weak supersolution $\phi$, just like in \cite{frank2008non} or (\ref{frank}). This allows us to obtain not only sharp remainders of weighted $L^{p}$-Hardy type inequalities, but also of $L^{p}$-Poincar\'e inequalities \cite[Corollary 3.5]{apseit2025sharp}: let $1<p<\infty$ and $\lambda_1$, $\phi_1$ be the first eigenvalue and eigenfunction of $-\Delta_{\gamma,p}$. Then,
\begin{align*}
\int_{\Omega}|\nabla_{\gamma}f|^{p}dz=\lambda_1\int_{\Omega}|f|^{p}dz+\int_{\Omega}C_p\left(\nabla_\gamma f,\phi_{1}\nabla_{\gamma}\left(\frac{f}{\phi_{1}}\right)\right)dz,
\end{align*}
for all complex-valued $f\in C^{\infty}_{0}(\Omega)$. Here $\Omega$ is a bounded open subset of $\mathbb{R}^{m+k}$ where the first eigenfunction $\phi_1$ of $-\Delta_{\gamma,p}$ is positive.

The condition (\ref{intro main condition}), in special cases (i.e. for radial weights), makes it possible to introduce the definition of Baouendi-Grushin $p$-Bessel pairs, which recovers the $p$-Bessel pairs definition introduced by Duy, Lam and Lu \cite{duy2022p}. We also refer to \cite{ruzhansky2024hardy} for related results. Last but not least, we note that all of our results are new in both Baouendi-Grushin and Euclidean frameworks.

The second aim of this paper is to apply the main result to obtain the sharp remainder formula of the Heisenberg-Pauli-Weyl inequality (in $L^{p}$) and consequently derive some stability results. First, let us recall the original extraordinary paper by Heisenberg \cite{heisenberg1927anschaulichen}, where he initially formulated the uncertainty principle. Essentially, Heisenberg discovered that the position and momentum of a particle cannot be measured at the same time. The idea was revolutionary and naturally led to its rigorous mathematical formulation by Kennard \cite{kennard1927quantenmechanik} and Weyl \cite{weyl1928gruppentheorie, weyl1950theory} (who credited it to Pauli). Nowadays, the following result is referred to as the Heisenberg–Pauli–Weyl (HPW) inequality; for consistency, we state it for compactly supported functions: let $f\in C^{\infty}_{0}(\mathbb{R}^{n})$. Then,
\begin{align}\label{hpw classical}
\left(\int_{\mathbb{R}^{n}}|\nabla f|^{2}dz\right)\left(\int_{\mathbb{R}^{n}}|z|^{2}|f|^{2}dz\right)\geq\frac{n^{2}}{4}\left(\int_{\mathbb{R}^{n}}|f|^{2}dz\right)^{2},
\end{align}
where the constant $\frac{n^2}{4}$ is optimal. Now the set of all extremizers of (\ref{hpw classical}) generally consists of Gaussian profiles of the form $f = ce^{-\alpha|z|^{2}}$ for some $c\in\mathbb{R}$ and $\alpha>0$. However, they are not compactly supported. So, in this context, they should actually be called \say{virtual} extremizers. However, if they describe all possible extremizers, we will still refer them as extremizers without saying that they are \say{virtual}. 

A natural question to ask is whether (\ref{hpw classical}) is stable or not, i.e. if there exists some \( f \) for which (\ref{hpw classical}) is close to an equality, does that imply that \( f \) is close to the set of extremizers? Recently, there has been a lot of progress made in this analysis. McCurdy and Venkatraman, in \cite{mccurdy2020quantitative}, were the first to address this question. The answer turns out to be positive. However, the paper \cite{mccurdy2020quantitative} had some missing gaps, which were filled in the subsequent papers \cite{fathi2021short, cazacu2024caffarelli, lam2025stability}. We also refer to \cite{cazacu2022sharp} for related results.

In the investigation of stability, one would like to have the explicit form of the so-called \say{Heisenberg deficit} $\delta_{2}(f)$:
\begin{align*}
\delta_{2}(f):=\left(\int_{\mathbb{R}^{n}}|\nabla f|^{2}dz\right)\left(\int_{\mathbb{R}^{n}}|z|^{2}|f|^{2}dz\right)-\frac{n^{2}}{4}\left(\int_{\mathbb{R}^{n}}|f|^{2}dz\right)^{2}.
\end{align*}
For example, in \cite[Corollary 2.6]{cazacu2024caffarelli}, Cazacu, Flynn, Lam and Lu obtained the following sharp remainder formula of a slightly different version of the HPW inequality:
\begin{align}\label{lu hpw}
\left(\int_{\mathbb{R}^{n}}|\nabla f|^{2}dz\right)^{\frac{1}{2}}\left(\int_{\mathbb{R}^{n}}|z|^{2}|f|^{2}dz\right)^{\frac{1}{2}}-\frac{n}{2}\int_{\mathbb{R}^{n}}|f|^{2}dz=\frac{\alpha^{2}}{2}\int_{\mathbb{R}^{n}}\left|e^{-\frac{|z|^{2}}{2\alpha^{2}}}\nabla\left(\frac{f}{e^{-\frac{|z|^{2}}{2\alpha^{2}}}}\right)\right|^{2}dz,
\end{align}
where $\alpha=\left(\frac{\int_{\mathbb{R}^{n}}|f|^{2}|z|^{2}dz}{\int_{\mathbb{R}^{n}}|\nabla f|^{2}dz}\right)^{\frac{1}{4}}$. The HPW identity (\ref{lu hpw}) was shown to provide valuable insights about the HPW inequality (\ref{hpw classical}). For example, the remainder term, in (\ref{lu hpw}), vanishes if and only if $f=ce^{-\frac{|z|^{2}}{2\alpha^{2}}}$ for any $c\in\mathbb{R}$. In fact, they proved that this implies that all extremizers of the classical HPW inequality (\ref{hpw classical}) are in the set $E=\{ce^{-\beta|z|^{2}}:c\in\mathbb{R},\beta>0\}$. Additionally, as shown in \cite{cazacu2024caffarelli}, the identity (\ref{lu hpw}) along with the optimal Poincar\'e inequality of Gaussian type measure (from \cite{do2023p})
\begin{align*}
\int_{\mathbb{R}^{n}}|\nabla f|^{2}e^{-\frac{|z|^{2}}{2|\alpha|^{2}}}dz\geq\frac{1}{|\alpha|^{2}}\inf_{c}\int_{\mathbb{R}^{n}}|f-c|^{2}e^{-\frac{|z|^{2}}{2|\alpha|^{2}}}dz
\end{align*}
gives the following sharp stability result of the altered HPW inequality:
\begin{align}\label{stability ineq}
\left(\int_{\mathbb{R}^{n}}|\nabla f|^{2}dz\right)^{\frac{1}{2}}\left(\int_{\mathbb{R}^{n}}|z|^{2}|f|^{2}dz\right)^{\frac{1}{2}}-\frac{n}{2}\int_{\mathbb{R}^{n}}|f|^{2}dz\geq d(f,E)^{2},
\end{align}
where $d(f,E)^{2}$ is the $L^{2}$ distance function from $f$ to the set $E$. Moreover, the inequality (\ref{stability ineq}) is sharp and the equality can be obtained by nontrivial functions $f\notin E$. In \cite{cazacu2024caffarelli}, the authors demonstrated that the stability properties of the altered HPW inequality actually imply and improve the corresponding stability properties of the classical HPW inequality (\ref{hpw classical}).

In this paper, we establish the sharp remainder formula and stability of the $L^{p}$-HPW inequality. Here, due to the novelty of the results, we will state them in the standard Euclidean framework (see Section \ref{sec hpw sharp} for its generalization to the Baouendi-Grushin setting, and \cite[Section 3.1.1]{KY18} for related results): let $1<p<\infty$ and $p'=\frac{p}{p-1}$. Then, for all non-zero complex-valued $f\in C^{\infty}_{0}(\mathbb{R}^{n})$, we have
\begin{multline}\label{p hpw identity}
\left(\int_{\mathbb{R}^{n}}|\nabla f|^{p}dz\right)\left(\int_{\mathbb{R}^{n}}|z|^{p'}|f|^{p}dz\right)^{p-1}-\left(\frac{n}{p}\right)^{p}\left(\int_{\mathbb{R}^{n}}|f|^{p}dz\right)^{p}
\\=\left(\int_{\mathbb{R}^{n}}|z|^{p'}|f|^{p}dz\right)^{p-1}\int_{\mathbb{R}^{n}}C_{p}\left(\nabla f,e^{-\alpha_{p,n}|z|^{p'}}\nabla\left(\frac{f}{e^{-\alpha_{p,n}|z|^{p'}}}\right)\right)dz,
\end{multline}
where $\alpha_{p,n}=\frac{n}{p}\frac{p-1}{p}\frac{\int_{\mathbb{R}^{n}}|f|^{p}dz}{\int_{\mathbb{R}^{n}}|z|^{p'}|f|^{p}dz}$. 

As in (\ref{lu hpw}), the identity (\ref{p hpw identity}) automatically gives the fact that all extremizers of the corresponding $L^{p}$-HPW inequality are of the form $f=ce^{-\alpha_{p,n}|z|^{p'}}$ for any $c\in\mathbb{C}$. We refer to \cite[Step 3 of Proof of Theorem 1.2]{cazacu2024hardy} and \cite[Lemmata 2.2 and 2.3]{CT24} for the estimates of the $C_p$-functional that allow us to make such claim. On top of that, since for any $f=ce^{\beta|z|^{p'}}$ ($\beta>0$) the formula for $\alpha_{p,n}$ gives $\alpha_{p,n}=\beta$, we have that the manifold of extremizers of the $L^{p}$-HPW inequality is $E_{p}=\{ce^{-\beta|z|^{p'}} : c\in\mathbb{C},\beta>0\}$.

When $p=2$, in (\ref{p hpw identity}), we derive the following sharp remainder of the $L^{2}$-HPW inequality (\ref{hpw classical}): for all non-zero complex-valued $f\in C^{\infty}_{0}(\mathbb{R}^{n})$
\begin{multline*}
\left(\int_{\mathbb{R}^{n}}|\nabla f|^{2}dz\right)\left(\int_{\mathbb{R}^{n}}|z|^{2}|f|^{2}dz\right)-\frac{n^2}{4}\left(\int_{\mathbb{R}^{n}}|f|^{2}dz\right)^{2}
\\=\left(\int_{\mathbb{R}^{n}}|z|^{2}|f|^{2}dz\right)\int_{\mathbb{R}^{n}}C_{2}\left(\nabla f,e^{-\alpha_{2,n}|z|^{2}}\nabla\left(\frac{f}{e^{-\alpha_{2,n}|z|^{2}}}\right)\right)dz
\\=\left(\int_{\mathbb{R}^{n}}|z|^{2}|f|^{2}dz\right)\int_{\mathbb{R}^{n}}\left|e^{-\alpha_{2,n}|z|^{2}}\nabla\left(\frac{f}{e^{-\alpha_{2,n}|z|^{2}}}\right)\right|^{2}dz,
\end{multline*}
where $\alpha_{2,n}=\frac{n}{4}\frac{\int_{\mathbb{R}^{n}}|f|^{2}dz}{\int_{\mathbb{R}^{n}}|z|^{2}|f|^{2}dz}$. Similarly, here we have $E=\{ce^{-\beta|z|^{2}} : c\in\mathbb{C}, \beta>0\}$.

Similarly as in \cite{cazacu2024caffarelli}, our $L^{p}$-HPW identity (\ref{p hpw identity}) along with weighted $L^{p}$-Poincar\'e inequality for the log-concave probability measure \cite[Lemma 4.1]{do2023p} gives us a quantitative stability result of the $L^{p}$-HPW inequality. First, let us define the $L^{p}$-Heisenberg deficit term:
\begin{align*}
\delta_{p}(f):=\left(\int_{\mathbb{R}^{n}}|\nabla f|^{p}dz\right)\left(\int_{\mathbb{R}^{n}}|z|^{p'}|f|^{p}dz\right)^{p-1}-\left(\frac{n}{p}\right)^{p}\left(\int_{\mathbb{R}^{n}}|f|^{p}dz\right)^{p}.
\end{align*}
Let $2\leq p<n$ and $p'=\frac{p}{p-1}$. Then, for all non-zero complex-valued $f\in C^{\infty}_{0}(\mathbb{R}^n)$, we have
\begin{align*}
\delta_{p}(f)\geq \widetilde{C}(n,p)\left(\int_{\mathbb{R}^{n}}|z|^{p'}|f|^{p}dz\right)^{p-1}d(f,E_{p})^{p},
\end{align*}
where $\widetilde{C}(n,p)>0$ and $d(f,E_{p})$ is the $L^{p}$ distance to the manifold of extremizers, i.e.
\begin{align*}
d(f,E_{p})^p:=\inf_{c\in\mathbb{C},\beta>0}\left\{\norm{f-ce^{-\beta|z|^{p'}}}^{p}_{L^{p}(\mathbb{R}^n)}\right\}.
\end{align*}

The final part of this paper will be concerned with the non-existence of positive solutions to the following nonlinear partial differential equation:
\begin{align}\label{PDE}
   \begin{cases}
\frac{\partial u}{\partial t} = \Delta^{v(z)}_{\gamma,p} u + \tilde{w}(z) u^{p-1} + \lambda u^q & \text{in } \Omega \times (0, T), \\
u(z, 0) = u_0(z) \geq 0 & \text{in } \Omega, \\
u(z, t) = 0 & \text{on } \partial \Omega \times (0, T),
\end{cases}
\end{align}
where $\Delta^{v(z)}_{\gamma,p} u=\nabla_\gamma \cdot\left(v(z)\left|\nabla_\gamma u \right|^{p-2} \nabla_\gamma u\right)$ is the weighted $p$-Grushin operator, $\lambda\in \mathbb{R}$, $q>0$ and $\tilde{w}(z)=\beta w(z)$ for $\beta>1$.

The motivation behind such PDE initially comes from the heat equation with singular potentials studied by Baras and Goldstein \cite{baras1984heat}:
\begin{align}\label{pde bg}
\begin{cases}
u_{t}=\Delta u+\frac{c}{|z|^{2}} u+f(z, t) & \text { if } z \in \Omega \subseteq \mathbb{R}^{n}, n \geq 3, t>0, c \in \mathbb{R}, f \geq 0, \\
u(z, 0)=u_{0}(z) \geq 0 & \text { if } x \in \Omega, u_{0} \in L^{2}(\Omega), \\
u(z, t)=0 & \text { if } z \in \partial \Omega, t>0,
\end{cases}    
\end{align}
where $\Omega$ is the whole space $\mathbb{R}^{n}$ or some bounded domain containing \\ $B_{1}=\left\{z \in \mathbb{R}^{n}:|z|<1\right\}$. There, they proved that if $c\leq \left(\frac{n-2}{2}\right)^{2}$, then (\ref{pde bg}) has a unique positive solution in the distributional sense if 
\begin{align*}
\int_{\Omega}|z|^{-\alpha} u_{0}(z) dz < \infty, \quad \int_{0}^{T} \int_{\Omega}|z|^{-\alpha} f(z, t) dzdt < \infty,
\end{align*}
where $\alpha$ is the smallest root of $\alpha^{2}-(N-2) \alpha+c=0$. Otherwise, when $c$ is greater than the sharp Hardy constant, then (\ref{pde bg}) has no positive solutions for $u_{0}>0$ and $f>0$. 

After that many authors have worked in developing the results of Baras and Goldstein \cite{baras1984heat} (see, \cite{azorero1998hardy, cabre1999existence, crespo2000global, kombe2004linear, goldstein2003linear, goldstein2003nonlinear, goldstein2005nonlinear, Kom06}). For instance, in \cite{crespo2000global}, Crespo and Alonso considered an $L^{p}$ analouge of (\ref{PDE}):
\begin{align*}
\begin{cases} 
\frac{\partial u}{\partial t} = \Delta_p u + \frac{c}{|z|^p} u^{p-1}, & \text{in } \mathbb{R}^n \times (0, \infty), \\
u(z, 0) = u_0(z) \geq 0 & \text{in } \mathbb{R}^n .
\end{cases}
\end{align*}
To be more precise, they showed that for the corresponding result to hold in $L^{p}$, one needs not only to control not only the constant $c$, but also the exponent $p$. In this direction, it is important to mention the earlier classical paper by Fujita \cite{fujita1966blowing}, where the impact of the exponent on the existence and non-existence of positive solution appears.

At the same time, Cabr\'e and Martel \cite{cabre1999existence} showed that the existence and non-existence of positive solutions to
\begin{align*}
\begin{cases}
\frac{\partial u}{\partial t}=\Delta u+V(z) u & \text { in } \quad \Omega \times(0, T), \\
u(z, 0)=u_{0}(z) \geq 0 & \text { in } \quad \Omega, \\
u(z, t)=0 & \text { on } \quad \partial \Omega \times(0, T)
\end{cases}
\end{align*}
is determined by the value of the infimum
\begin{align*}
\inf_{0 \neq \phi \in C_{0}^{\infty}(\Omega)} \frac{\int_{\Omega}|\nabla \phi|^{2} dz - \int_{\Omega} V|\phi|^{2} dz}{\int_{\Omega}|\phi|^{2} dz}.
\end{align*}

The technique was later shown to be extremely versatile in a vast range of different settings, such as \cite{goldstein2003instantaneous, kombe2004linear, goldstein2005nonlinear, Kom06, goldstein2022nonexistence} just to name a few. In one of them \cite{Kom06}, by combining the ideas from \cite{crespo2000global} and \cite{cabre1999existence}, Kombe derived a similar condition of Cabr\'e and Martel that implied the non-existence of positive solutions to the nonlinear degenerate parabolic equation:
\begin{align}\label{kombe pde}
\begin{cases} 
\frac{\partial u}{\partial t} = \nabla_{\gamma} \cdot (|\nabla_{\gamma} u|^{p-2} \nabla_{\gamma} u) + V(z) u^{p-1} & \text{in} \quad \Omega \times (0,T), \quad 1 < p < 2, \\ 
u(z,0) = u_0(z) \geq 0 & \text{in} \quad \Omega, \\ 
u(z,t) = 0 & \text{on} \quad \partial \Omega \times (0,T)
\end{cases}
\end{align}
in the framework of Baouendi-Grushin vector fields. Then, in \cite{han2015nonlinear}, Han and Guo provided several non-existence results to equations of the type (\ref{kombe pde}) with weights in the divergence term $\nabla_{\gamma}\cdot\left(\rho^{-p\lambda}|\nabla_{\gamma}u|^{p-2}\nabla_{\gamma}u\right)$ and time-dependent singular potentials $V(z,t)$.

In this paper, we study the non-existence of positive solutions to the PDE (\ref{PDE}). We note that our approach is reminiscent of that given in \cite{Kom06} and \cite{goldstein2022nonexistence}. We then apply the general condition (Theorem \ref{parabolic eq}) to study the non-existence of positive solution to specific PDEs (Corollaries \ref{cor1} and \ref{cor2}). 

The paper is organized as follows. Section~\ref{sec main results} presents our main results and their consequences for general weights. In Section~\ref{sec applications}, we apply Theorem~\ref{main thm} to recover the sharp remainders of most known \(L^{p}\)-Hardy-type inequalities (Section \ref{sec sharp remainders}), covering both radial and non-radial weights. We also establish new identities for non-radial weights and derive the sharp remainder term for the \(L^{p}\)-Poincar\'e inequality. Section \ref{sec hpw sharp} shows the sharp remainder and stability of the $L^{p}$-HPW inequality. Finally, in Section \ref{sec: non-exist}, we derive the condition (Theorem \ref{parabolic eq}), with two Corollaries (Corollary \ref{cor1} and Corollary \ref{cor2}), under which the non-existence of positive local solutions to (\ref{PDE}) holds. Section~\ref{proof of main} contains the proof of the main result, and Section~\ref{proofs of cor} provides the proofs of all applications of the main result.

\section{Main results}\label{sec main results}
In this section, we present the main results of the paper. Before stating the main results, we recall the following definition and estimates of the $C_p$-functional that will be used later. 

\begin{defn}
Let $1<p<\infty$. Then, for $\xi,\eta\in\mathbb{C}^{n}$, we define
\begin{align*}
C_p(\xi,\eta):=|\xi|^p-|\xi-\eta|^p-p|\xi-\eta|^{p-2}\textnormal{Re}(\xi-\eta)\cdot\overline{\eta}\geq0.
\end{align*}
\end{defn}

\begin{lem}[\text{\cite[Step 3 of Proof of Theorem 1.2]{cazacu2024hardy}}]\label{lem1}
Let $p\geq2$. Then, for $\xi,\eta\in\mathbb{C}^n$, we have
\begin{align*}
C_{p}(\xi,\eta)\geq c_1(p)|\eta|^{p},
\end{align*}
where 
\begin{align}\label{c1 const}
c_1(p)
= \inf_{(s,t)\in\mathbb{R}^2\setminus\{(0,0)\}}
\frac{\bigl[t^2 + s^2 + 2s + 1\bigr]^{\frac p2} -1-ps}
{\bigl[t^2 + s^2\bigr]^{\frac p2}}\in(0,1].
\end{align}
\end{lem}
\begin{lem}[\text{\cite[Lemma 2.2]{CT24}}]\label{lem2}
Let $1<p<2\leq n$. Then, for $\xi,\eta\in \mathbb{C}^{n}$, we have
\begin{align*}
C_p(\xi, \eta) \geq c_2(p) \frac{|\eta|^2}{\left( |\xi| + |\xi - \eta| \right)^{2-p}},
\end{align*}
where
\begin{align}\label{c2 const}
c_2(p) := \sup_{s^2 + t^2 > 0} \frac{\left( t^2 + s^2 + 2s + 1 \right)^{\frac{p}{2}} - 1 - ps}{\left( \sqrt{t^2 + s^2 + 2s + 1} + 1 \right)^{p-2} (t^2 + s^2)} \in \left(0,  \frac{p(p-1)}{2^{p-1}} \right].
\end{align}
\end{lem}
\begin{lem}[\text{\cite[Lemma 2.3]{CT24}}]\label{lem3}
Let $1<p<2\leq n$. Then, for $\xi,\eta\in \mathbb{C}^{n}$, we have
\begin{align*}
C_p(\xi, \eta) \leq c_3(p) \frac{|\eta|^2}{\left( |\xi| + |\xi - \eta| \right)^{2-p}},
\end{align*}
where
\begin{align}\label{c3 const}
c_3(p) := \sup_{s^2 + t^2 > 0} \frac{\left( t^2 + s^2 + 2s + 1 \right)^{\frac{p}{2}} - 1 - ps}{\left( \sqrt{t^2 + s^2 + 2s + 1} + 1 \right)^{p-2} (t^2 + s^2)} \in \left[ \frac{p}{2^{p-1}}, +\infty \right).
\end{align}
\end{lem}

We now state the main result of this section. 

\begin{thm}\label{main thm}
Let $1<p<\infty$ and let $\Omega\subseteq\mathbb{R}^{m+k}$ be an open set such that the integrals below make sense. Let $v,w\in L^{1}_{loc}(\Omega)$ be nonnegative functions with $\phi>0$ satisfying the following in the distributional sense in $\Omega$:
\begin{align}\label{main condition}
-\nabla_{\gamma} \cdot \left( v \left|\nabla_{\gamma} \phi\right|^{p-2} \nabla_{\gamma} \phi \right)\geq w\phi^{p-1}.
\end{align}
Then, for all complex-valued functions $f\in C_{0}^{\infty}(\Omega)$, we have, 
\begin{align}\label{main result}
\int_{\Omega}^{}v|\nabla_{\gamma}f|^{p}dz \geq \int_{\Omega}^{}w|f|^{p}dz + \int_{\Omega}^{}vC_{p}\left( \nabla_{\gamma}f,\phi\nabla_{\gamma}\left( \frac{f}{\phi} \right)  \right) dz,
\end{align}
Equality, in (\ref{main result}), holds when $\phi$ is a positive distributional solution to (\ref{main condition}) with an equality. 
\end{thm}

\begin{comment}

\begin{rem}
In this paper, we shall mainly focus on the case when $-\nabla_{\gamma} \cdot \bigl( v \left|\nabla_{\gamma} \phi\right|^{p-2}$ $\nabla_{\gamma} \phi \bigr)$ $=$ $w\phi^{p-1}$. In particular, given any weighted functional inequality of the following form
\begin{align}\label{hardy weighted}
\int_{\Omega}v|\nabla_{\gamma}f|^{p}dz\geq\int_{\Omega}w|f|^{p}dz
\end{align}
we will take suitable $v$ and $\phi$ such that after calculating 
\begin{align}\label{main condition =}
w=\frac{-\nabla_{\gamma}\cdot \left( v|\nabla_{\gamma}\phi|^{p-2}\nabla_{\gamma}\phi \right) }{\phi^{p-1}}
\end{align}
we get $w$ from (\ref{hardy weighted}) with (sometimes) additional remainder terms. This approach allows us to obtain identities with positive remainder terms for the standard sub-elliptic Grushin operator $\nabla_{\gamma}$.
\end{rem}

\end{comment}

When $\gamma=0$, in Theorem \ref{main thm}, the result is also new in the Euclidean setting.

\begin{cor}%\label{cor of thm}
Let $1<p<\infty$ and let $\Omega\subseteq\mathbb{R}^{n}$ be an open set such that the integrals below make sense. Let $v,w\in L^{1}_{loc}(\Omega)$ be nonnegative functions with $\phi>0$ satisfying the following in the distributional sense in $\Omega$:
\begin{align}\label{euc main condition}
-\nabla \cdot \left( v \left|\nabla \phi\right|^{p-2} \nabla \phi \right)\geq w\phi^{p-1}.
\end{align}
Then, for all complex-valued functions $f\in C_{0}^{\infty}(\Omega)$, we have, 
\begin{align}\label{euc main result}
\int_{\Omega}^{}v|\nabla f|^{p}dz \geq \int_{\Omega}^{}w|f|^{p}dz + \int_{\Omega}^{}vC_{p}\left( \nabla f,\phi\nabla \left( \frac{f}{\phi} \right)  \right) dz.
\end{align}
Equality, in (\ref{euc main result}), holds when $\phi$ is a positive distributional solution to (\ref{euc main condition}) with an equality.
\end{cor}

The inequality (\ref{main result}) not only improves the results of Kombe and Yener (\ref{kom1}) and (\ref{kom2}), but also \say{unifies} the two weighted Hardy type inequalities, in the sense that it contains the two inequalities, (\ref{kom1}) and (\ref{kom2}), at once. Moreover, the inequalities are with explicit constants.

\begin{cor}\label{cor rec KY}
Let $\Omega$, $v, w$ and $\phi$ be from Theorem \ref{main thm}. Then, for $p\geq2$, we have 
\begin{align*}
\int_{\Omega} v\left|\nabla_\gamma f\right|^p dz \geq \int_{\Omega} w|f|^p dz+c_1(p) \int_{\Omega} v\left|\nabla_\gamma \left(\frac{f}{\phi}\right)\right|^p \phi^p dz,
\end{align*}
where the constant $c_1(p)$ is given in (\ref{c1 const}). In addition, for $1<p<2\leq n$, we have
\begin{align*}
\int_{\Omega}v\left|\nabla_\gamma f\right|^p dz \geq \int_{\Omega} w|f|^p dz+\widetilde{c_2}(p) \int_{\Omega} \frac{v\left|\nabla_\gamma \left(\frac{f}{\phi}\right)\right|^2 \phi^2}{\left(\left|\frac{f}{\phi} \nabla_\gamma \phi\right|+\left|\nabla_\gamma \left(\frac{f}{\phi}\right)\right| \phi\right)^{2-p}} dz,
\end{align*}
where $\widetilde{c_2}(p)=2^{p-2}c_2(p)$ with $c_2(p)$ given in (\ref{c2 const}).
\end{cor}

\begin{rem}
If $v, w$ and $\phi$ are radial (i.e. depend only on $\rho$), then the condition (\ref{main condition}) can be reduced to the $p$-Bessel pair (in the setting of Baouendi-Grushin vector fields) first introduced in \cite{duy2022p} in the standard Euclidean setting. Here, we show this implication: suppose we have
\begin{align}\label{bess1}
-\nabla_{\gamma} \cdot \left( v \left|\nabla_{\gamma} \phi\right|^{p-2} \nabla_{\gamma} \phi \right) = w \phi^{p-1}
\end{align}
with $v=v(\rho)$, $w=w(\rho)$ and $\phi=\phi(\rho)$. Then, we have
\begin{align}\label{bess2}
\nabla_{\gamma}\phi=\phi'(\rho)\nabla_{\gamma}\rho \implies |\nabla_{\gamma}\phi|=|\phi'(\rho)||\nabla_{\gamma}\rho|. 
\end{align}
Using (\ref{bess2}), we can rewrite (\ref{bess1}) as
\begin{align}\label{bess3}
-\nabla_{\gamma}\cdot\left(v(\rho)|\phi'(\rho)|^{p-2}\phi'(\rho)|\nabla_{\gamma}\rho|^{p-2}\nabla_{\gamma}\rho\right)=w(\rho)\phi^{p-1}.
\end{align}
To simplify the calculations, let us denote
\begin{align*}
A(\rho):=v|\phi'(\rho)|^{p-2}\phi'(\rho).
\end{align*}
Now we calculate
\begin{align}\label{bess4}
\nabla_{\gamma}\cdot\left(A(\rho)|\nabla_{\gamma}\rho|^{p-2}\nabla_{\gamma}\rho\right)&=A'(\rho)|\nabla_{\gamma}\rho|^{p}+A(\rho)\Delta_{\gamma,p}\rho \nonumber
\\&=A'(\rho)|\nabla_{\gamma}\rho|^{p}+\frac{Q-1}{\rho}A(\rho)|\nabla_{\gamma}\rho|^{p} \nonumber
\\&=\frac{|\nabla_{\gamma}\rho|^p}{\rho^{Q-1}}\left(\rho^{Q-1} A(\rho)\right)'.
\end{align}
Substituting (\ref{bess4}) into (\ref{bess3}), we get
\begin{multline*}
-\nabla_{\gamma}\cdot\left(v(\rho)|\phi'(\rho)|^{p-2}\phi'(\rho)|\nabla_{\gamma}\rho|^{p-2}\nabla_{\gamma}\rho\right)=-\frac{|\nabla_{\gamma}\rho|^p}{\rho^{Q-1}}\left(\rho^{Q-1} v|\phi'(\rho)|^{p-2}\phi'(\rho)\right)'\\=w(\rho)\phi^{p-1},
\end{multline*}
which gives us
\begin{align*}%\label{final bessel}
\left(\rho^{Q-1} v|\phi'(\rho)|^{p-2}\phi'(\rho)\right)'+\frac{\rho^{Q-1}}{|\nabla_{\gamma}\rho|^p}w(\rho)\phi^{p-1}=0.
\end{align*}
If we do not assume that $\phi$ is positive, then (\ref{bess1}) becomes
\begin{align*}%\label{rec1}
-\nabla_{\gamma} \cdot \left( v \left|\nabla_{\gamma} \phi\right|^{p-2} \nabla_{\gamma} \phi \right) = w |\phi|^{p-2}\phi.
\end{align*}
By the same process, we obtain
\begin{align}\label{final final bessel}
\left(\rho^{Q-1} v(\rho)|\phi'(\rho)|^{p-2}\phi'(\rho)\right)'+\frac{\rho^{Q-1}}{|\nabla_{\gamma}\rho|^p}w(\rho)|\phi|^{p-2}\phi=0.
\end{align}
If $\gamma=0$ in (\ref{final final bessel}), we recover the $p$-Bessel pair definition from \cite[Definition 1.1]{duy2022p} by Duy, Lam and Lu. Therefore, we can also introduce the following definition of $p$-Bessel pairs, now, in the setting of Baouendi-Grushin vector fields.
\end{rem}

\begin{defn}\label{p-bessel pair}
The pair $(v(\rho),w(\rho))$ is called a Baouendi-Grushin $p$-Bessel pair on $(0,R)$ if  
\[
\left(\rho^{Q-1} v(\rho)|\phi'(\rho)|^{p-2}\phi'(\rho)\right)' + \frac{\rho^{Q-1}}{|\nabla_{\gamma}\rho|^p}w(\rho)|\phi|^{p-2}\phi = 0
\]
has a positive solution on $(0,R)$.
\end{defn}

\section{Applications of Theorem \ref{main thm}}\label{sec applications}

In this section, we apply the Theorem \ref{main thm}. However, before we start applying the Theorem \ref{main thm}, we note that to make the proofs rigorous, one should (for any given $\epsilon>0$) replace the function $\rho$ by its regularization
\begin{align*}
\rho_{\epsilon}:=\left(|x|^{2(1+\gamma)}_{\epsilon}+(1+\gamma)^{2}|y|^{2}\right)^{\frac{1}{2(1+\gamma)}}
\end{align*}
with $|x|_{\epsilon}:=\left(\epsilon^{2}+\sum_{i=1}^{m}x^{2}_{i}\right)^{\frac{1}{2}}$ and then, after calculations take limit $\epsilon \rightarrow 0$. However, for the simplicity, we proceed formally.

Lastly, to clarify, we note that the results below provide new identities to known inequalities. When the Baouendi-Grushin vector fields are reduced to the standard Euclidean gradient, all of the results below are also new.

\subsection{Sharp remainders of weighted Hardy-type, Poincar\'e and Hardy-Poin-car\'e inequalities}\label{sec sharp remainders} Now as previously mentioned, we can employ Theorem \ref{main thm} to obtain sharp remainder terms of known weighted Hardy type, Poincar\'e and Hardy-Poincar\'e inequalities associated to the Baouendi-Grushin operator. We begin with
\begin{align*}
v = 1 \quad \text{and} \quad \phi=(R-\rho)^{\frac{p-1}{p}},
\end{align*}
which gives us the sharp remainder terms of the Hardy type inequality obtained by Niu, Chen and Han \cite[Theorem 2.1]{NCH04}:

\begin{cor}\label{cor nhc}
Let $1<p<\infty$ and $R>0$. Then, for all complex-valued $f\in C^{\infty}_{0}(B_{R}\backslash\{(0,0)\})$, we have
\begin{multline*}
\int_{B_{R}}^{}|\nabla_{\gamma}f|^{p}dz=\left( \frac{p-1}{p} \right) ^{p}\int_{B_{R}}^{}\frac{|x|^{\gamma p}}{(R-\rho)^{p}\rho^{\gamma p}}|f|^{p}dz\\+\int_{B_{R}}^{}C_{p}\left( \nabla_{\gamma}f,(R-\rho)^{\frac{p-1}{p}}\nabla_{\gamma}\left(  \frac{f}{(R-\rho)^{\frac{p-1}{p}}}\right)  \right)dz\\+(Q-1)\left( \frac{p-1}{p} \right) ^{p-1}\int_{B_{R}}^{}\frac{|x|^{\gamma p}}{(R-\rho)^{p-1}\rho^{\gamma p+1}}|f|^{p}dz.
\end{multline*}
\end{cor}

We can also obtain the sharp remainder term of the result by D'Ambrosio \cite[Formula (3.1) of Theorem 3.1]{D'A04}. In particular, if we set
\begin{align*}
v=|x|^{\beta-\gamma p}\rho^{(1+\gamma)p-\alpha} \quad \text{and} \quad \phi=\rho^{-\frac{Q+\beta -\alpha}{p}},
\end{align*}
then we get $w=\left( \frac{Q+\beta-\alpha}{p} \right)^{p}\frac{|x|^{\beta}}{\rho^{\alpha}}$, and therefore the following result:

\begin{cor}\label{cor dam}
Let $1<p<\infty$, $m,k\geq1$ and $\alpha,\beta \in\mathbb{R}$ such that $Q\geq\alpha-\beta$. Then, for all complex-valued $f\in C^{\infty}_{0}(\mathbb{R}^{m+k}\backslash\{(0,0)\})$, we have
\begin{multline}\label{dam org eq}
\int_{\mathbb{R}^n}^{}|x|^{\beta-\gamma p}\rho^{(1+\gamma)p-\alpha}|\nabla_{\gamma}f|^{p}dz=\left( \frac{Q+\beta-\alpha}{p} \right) ^{p}\int_{\mathbb{R}^n}^{}\frac{|x|^{\beta}}{\rho^{\alpha}}|f|^{p}dz\\+\int_{\mathbb{R}^n}^{}|x|^{\beta-\gamma p}\rho^{(1+\gamma)p-\alpha}C_{p}\left( \nabla_{\gamma}f,\rho^{-\frac{Q+\beta-\alpha}{p}}\nabla_{\gamma}\left(  \frac{f}{\rho^{-\frac{Q+\beta-\alpha}{p}}}\right)  \right)dz.
\end{multline}
\end{cor}

\begin{rem}
Due to the importance of \cite[Formula (3.1) of Theorem 3.1]{D'A04}, a lot of consequent results emerge out of (\ref{dam org eq}). For instance, when $\alpha=(1+\gamma)p$ and $\beta=\gamma p$ in (\ref{dam org eq}), we can obtain the sharp remainder formula of \cite[Formula (3.2) of Theorem 3.1]{D'A04}. That is, we have
\begin{align}\label{new res}
\int_{\mathbb{R}^n}|\nabla_{\gamma}f|^{p}dz=\left(\frac{Q-p}{p}\right)^{p}\int_{\mathbb{R}^n}\frac{|x|^{\gamma p}}{\rho^{\gamma p+p}}|f|^{p}dz+\int_{\mathbb{R}^n}C_{p}\left(\nabla_{\gamma}f,\rho^{-\frac{Q-p}{p}}\nabla_{\gamma}\left(\frac{f}{\rho^{-\frac{Q-p}{p}}}\right)\right)dz,
\end{align}
which for $p=2$, in (\ref{new res}), gives the sharp remainder term of (\ref{Gar}). Furthermore, when $\gamma=0$, we recover the sharp remainder term of the classical $L^{p}$-Hardy inequality (\ref{class_multdim_hardy}) from \cite[Lemma 3.3]{cazacu2024hardy}.
\end{rem}

We also obtain a particular case of the result from \cite[Formula (3.11) of Theorem 3.5]{D'A05}. That is, picking
\begin{align*}
v=|x|^{\alpha+p} \quad \text{and} \quad \phi=|x|^{\frac{|m+\alpha|}{p}}
\end{align*}
gives exactly $w=\left(\frac{|m+\alpha|}{p}\right)^{p}|x|^{\alpha}$. Hence, the next corollary holds.

\begin{cor}\label{cor dam |x|}
Let $1<p<\infty$, $\alpha<-m$ and $\Omega\subset(\mathbb{R}^{m}\backslash\{0\})\times\mathbb{R}^k$ be an open set. Then, we have 
\begin{multline}\label{dam |x| eq}
\int_{\Omega} |x|^{\alpha+p}\left|\nabla_{\gamma}f\right|^p dz =\left(\frac{|m+\alpha|}{p}\right)^{p} \int_{\Omega}|x|^{\alpha}|f|^pdz
\\+\int_{\Omega} |x|^{\alpha+p}C_p\left(\nabla_{\gamma}f, |x|^{\frac{|m+\alpha|}{p}} \nabla_{\gamma}\left(\frac{f}{|x|^{\frac{|m+\alpha|}{p}}}\right)\right)dz
\end{multline}
for all complex-valued $f\in C^{\infty}_{0}(\Omega)$.
\end{cor}

Apart from Hardy type inequalities, we also obtain the sharp remainder formula of the $L^{p}$-Poincar\'e inequality from \cite[Corollary 3.5]{apseit2025sharp} by specifying
\begin{align*}
v=1 \quad \text{and} \quad \phi=\phi_1:=\text{the first eigenfunction of } -\Delta_{\gamma,p},
\end{align*}
we get precisely $w=\lambda_{1}$. In this context, $\lambda_{1}$ corresponds to the first eigenvalue of the minus Dirichlet $p$-Grushin operator $-\Delta_{\gamma,p}$ on the open bounded subset $\Omega$ of $\mathbb{R}^{m+k}$ where the first eigenfunction is positive.
\begin{cor}\label{cor poin}
Let $1<p<\infty$ and $\Omega\subset\mathbb{R}^{m+k}$ . Then,
\begin{align*}
\int_{\Omega}|\nabla_{\gamma}f|^{p}dz=\lambda_1\int_{\Omega}|f|^{p}dz+\int_{\Omega}C_p\left(\nabla_\gamma f,\phi_{1}\nabla_{\gamma}\left(\frac{f}{\phi_{1}}\right)\right)dz,
\end{align*}
for all complex-valued $f\in C^{\infty}_{0}(\Omega)$. Here, $\Omega$ is a bounded open subset of $\mathbb{R}^{m+k}$ where the first eigenfunction $\phi_1$ of $-\Delta_{\gamma,p}$ is positive.
\end{cor}

Returning to Hardy inequalities, we observe that choosing
\begin{align*}
v=\left(\log\frac{R}{\rho}\right)^{\alpha+p} \quad \text{and} \quad \phi=\left(\log \frac{R}{\rho}\right)^{\frac{|\alpha+1|}{p}}
\end{align*}
produces the Hardy weight and an extra remainder term. That is, we get 
\begin{align*}
w=\left(\frac{|\alpha+1|}{p}\right)^p\left(\log\frac{R}{\rho}\right)^{\alpha}\frac{|x|^{\gamma p}}{\rho^{\gamma p + p}}+\left(\frac{|\alpha+p|}{p}\right)^{p-1}(Q-p)\left(\log\frac{R}{\rho}\right)^{\alpha+1}\frac{|x|^{\gamma p}}{\rho^{\gamma p +p}}
\end{align*}
and, thus, the following sharp remainder formula of the power logarithmic $L^{p}$-Hardy type inequality by D'Ambrosio \cite[Formula (3.8) of Theorem 3.1]{D'A05}:
\begin{cor}\label{cor log rho}
Let $1<p<\infty$, $\alpha<-1$ and $Q\geq p$. Then, we have
\begin{multline*}
\int_{B_{R}}^{}\left(\log\frac{R}{\rho}\right)^{\alpha+p}|\nabla_{\gamma}f|^{p}dz = \left( \frac{|\alpha+1|}{p} \right) ^{p}\int_{B_{R}}^{}\left( \log \frac{R}{\rho} \right)^{\alpha}\frac{|x|^{\gamma p}}{\rho^{\gamma p+p}}|f|^{p}dz\\+\left( \frac{|\alpha+1|}{p} \right) ^{p-1}(Q-p)\int_{B_{R}}^{}\left( \log \frac{R}{\rho}\right)^{\alpha+1}\frac{|x|^{\gamma p}}{\rho^{\gamma p+p}}|f|^{p}dz
\\+\int_{B_{R}}^{} \left(\log\frac{R}{\rho}\right)^{\alpha+p}C_{p}\left( \nabla_{\gamma}f,\left( \log \frac{R}{\rho}\right) ^{\frac{|\alpha+1|}{p}}\nabla_{\gamma}\left( \frac{f}{\left( \log \frac{R}{\rho}\right) ^{\frac{|\alpha+1|}{p}}} \right)  \right)dz
\end{multline*}
for all complex-valued $f\in C^{\infty}_{0}(B_{R}\backslash\{(0,0)\})$.
\end{cor}

A similar identity can be obtained, but now with different logarithmic weights, which are from \cite[Formula (3.13) of Theorem 3.5]{D'A05}. So, if we take
\begin{align*}
v=\left( \log \frac{R}{|x|} \right)^{\alpha+p} \quad \text{and} \quad \phi=\left( \log \frac{R}{|x|} \right)^{\frac{|\alpha+1|}{p}},
\end{align*}
we get $w=\left(\frac{|\alpha+1|}{p}\right)^{p}\frac{\left(\log \frac{R}{|x|}\right)^{\alpha}}{|x|^{p}}
+\left(\frac{|\alpha+1|}{p}\right)^{p-1}(m-p)\frac{\left(\log \frac{R}{|x|}\right)^{\alpha+1}}{|x|^{p}}$.
\begin{cor}\label{cor log |x|}
Let $1<p<\infty$, $\alpha<-1$ and $m\geq p$. Then, for all complex-valued $f\in C^{\infty}_{0}(\Omega\backslash\{x=0\})$, we have
\begin{multline*}
\int_{\Omega} \left(\log\frac{R}{|x|}\right)^{\alpha+p} |\nabla_{\gamma} f|^{p} dz 
= \left( \frac{|\alpha+1|}{p} \right)^{p} \int_{\Omega} \left(\log\frac{R}{|x|}\right)^{\alpha} \frac{|f|^{p}}{|x|^{p}} dz \\ + \left( \frac{|\alpha+1|}{p} \right)^{p-1} (m-p) \int_{\Omega} \left(\log\frac{R}{|x|}\right)^{\alpha+1} \frac{|f|^{p}}{|x|^{p}} dz \\+ \int_{\Omega} \left(\log\frac{R}{|x|}\right)^{\alpha+p} 
C_{p} \left( \nabla_{\gamma} f, \left( \log\frac{R}{|x|} \right)^{\frac{|\alpha+1|}{p}} 
\nabla_{\gamma} \left( \frac{f}{\left( \log\frac{R}{|x|} \right)^{\frac{|\alpha+1|}{p}}} \right) \right) dz,
\end{multline*}
where $\Omega=\{z\in\mathbb{R}^{m+k}:|x|<R\}$.
\end{cor}

By choosing
\begin{align*}
v=\left( 1+\rho^{\frac{p}{p-1}} \right)^{\alpha(p-1)} \quad \text {and} \quad \phi=\left( 1+\rho^{\frac{p}{p-1}} \right)^{1-\alpha},
\end{align*}
we get $w=Q\left(\frac{(\alpha-1)p}{p-1}\right)^{p-1}\left(1+\rho^{\frac{p}{p-1}}\right)^{(\alpha-1)(p-1)}\frac{|x|^{\gamma p}}{\rho^{\gamma p}}$, which allows us to derive the sharp remainder formula of the inequality proved in \cite[Corollary 3.8]{KY18} by Kombe and Yener.

\begin{cor}\label{cor rho alpha}
Let $1<p<\infty$ and $\alpha>1$. Then, for all complex-valued $f\in C^{\infty}_{0}(\mathbb{R}^{m+k}\backslash\{(0,0)\})$, we have
\begin{multline}\label{rho alpha eq}
\int_{\mathbb{R}^{n}}\left(1+\rho^{\frac{p}{p-1}}\right)^{\alpha(p-1)}|\nabla_{\gamma}f|^{p}dz\\=Q\left(\frac{(\alpha-1)p}{p-1}\right)^{p-1}\int_{\mathbb{R}^{n}}\left(1+\rho^{\frac{p}{p-1}}\right)^{(\alpha-1)(p-1)}\frac{|x|^{\gamma p}}{\rho^{\gamma p}}|f|^{p}dz\\+\int_{\mathbb{R}^{n}}\left(1+\rho^{\frac{p}{p-1}}\right)^{\alpha(p-1)}C_{p}\left(\nabla_{\gamma}f,\left( 1+\rho^{\frac{p}{p-1}} \right)^{1-\alpha}\nabla_{\gamma}\left(\frac{f}{\left( 1+\rho^{\frac{p}{p-1}} \right)^{1-\alpha}}\right)\right)dz.
\end{multline}
\end{cor}

\begin{rem}
Under $\gamma=0$, in (\ref{rho alpha eq}), we obtain the sharp remainder formula of the Hardy-Poincar\'e inequalities by Skrzypczak from \cite[Theorem 5.3]{skrzypczak2013hardy} and \cite[Theorem 3.1]{skrzypczak2014hardy}. We note that the inequality form of (\ref{rho alpha eq}) has many applications in the theory of nonlinear diffusions \cite{carrillo2002poincare, blanchet2007hardy, agueh2010large, dolbeault2011fast}. Also, we refer the reader to \cite[Section 4.1]{skrzypczak2014hardy} for the comparisons with the classical Hardy inequality and to \cite[Section 4.2]{skrzypczak2014hardy} for the connections between (\ref{rho alpha eq}) and the results from \cite{blanchet2007hardy, bonforte2010sharp, ghoussoub2011bessel}.
\end{rem}

We also establish sharp remainder terms of Hardy-type inequalities by Yener \cite[Theorem 3.3]{yener2018several}. That is, considering
\begin{align*}
v=\frac{(a+b\rho^{\alpha})^{\beta}}{\rho^{\ell p}} \quad \text{and} \quad \phi=\rho^{-\frac{Q-p\ell-p}{p}}
\end{align*}
we derive that 
\begin{multline*}
w=\left(\frac{Q-p\ell-p}{p}\right)^{p}(a+b\rho^{\alpha})^{\beta}\frac{|x|^{\gamma p}}{\rho^{\gamma p+\ell p+p}}
\\+\left(\frac{Q-p\ell-p}{p}\right)^{p-1}\beta\alpha b(a+b\rho^{\alpha})^{\beta-1}\frac{|x|^{\gamma p}}{\rho^{\gamma p+\ell p+p-\alpha}},
\end{multline*}
which, as a result, gives:
\begin{cor}\label{cor super}
Let $a,b>0$ and $\alpha,\beta, \ell$ be real numbers such that $\alpha\beta>0$ and $Q\geq pl+p$. Then, for all complex-valued $f\in C^{\infty}_{0}(\mathbb{R}^{m+k}\backslash\{(0,0)\})$, we have
\begin{multline}\label{super eq}
\int_{\mathbb{R}^{n}}\frac{(a+b\rho^{\alpha})^{\beta}}{\rho^{\ell p}}
|\nabla_{\gamma}f|^{p}dz
=\left(\frac{Q-p\ell-p}{p}\right)^{p}\int_{\mathbb{R}^{n}}
\frac{(a+b\rho^{\alpha})^{\beta}}{\rho^{\ell p+p}}\frac{|x|^{\gamma p}}{\rho^{\gamma p}}|f|^{p}dz
\\+\left(\frac{Q-p\ell-p}{p}\right)^{p-1}\beta\alpha b\int_{\mathbb{R}^{n}}
\frac{(a+b\rho^{\alpha})^{\beta-1}}{\rho^{\ell p+p-\alpha}}\frac{|x|^{\gamma p}}{\rho^{\gamma p}}|f|^{p}dz
\\+\int_{\mathbb{R}^{n}}\frac{(a+b\rho^{\alpha})^{\beta}}{\rho^{\ell p}}
C_{p}\left(\nabla_{\gamma}f,\rho^{-\frac{Q-p\ell-p}{p}}\nabla_{\gamma}\left(\frac{f}{\rho^{-\frac{Q-p\ell-p}{p}}}\right)\right)dz.
\end{multline}
\end{cor}

\begin{rem}
In the $L^{2}$ Euclidean case ($\gamma=0$), the identity (\ref{super eq}) gives sharp remainder terms of the weighted Hardy inequalities by Ghoussoub and Moradifam \cite[Formula (7)]{ghoussoub2011bessel}. For the applications of (\ref{super eq}) in the theory of ODE, we refer to \cite[Theorem 2.1]{ghoussoub2011bessel}.
\end{rem}

In the spirit of \cite{yener2018several}, we also prove several sharp remainder formulas of $L^{p}$-Hardy inequalities related to the Baouendi-Grushin operator with other non-radial weights. For example, setting
\begin{align*}
v=\left(\frac{y_1}{|x|^{\gamma}}\right)^{p-2}\log x_{1} \quad \text{and} \quad \phi=\log y_{1},
\end{align*}
gives the following sharp remainder formula of \cite[Theorem 3.6]{yener2018several}:
\begin{cor}\label{cor several yener}
Let $1<p<\infty$. Then, for all complex-valued $f\in C^{\infty}_{0}(\Omega)$, we have
\begin{multline*}
\int_{\Omega}\left(\frac{y_1}{|x|^{\gamma}}\right)^{p-2}\log x_{1}|\nabla_{\gamma}f|^{p}dz=\int_{\Omega}\frac{|x|^{2\gamma}\log x_{1}}{y^{2}_{1}\log^{p-1}y_{1}}|f|^{p}dz
\\+\int_{\Omega}\left(\frac{y_1}{|x|^{\gamma}}\right)^{p-2}\log x_{1}C_{p}\left(\nabla_{\gamma}f,\log y_{1}\nabla_{\gamma}\left(\frac{f}{\log y_{1}}\right)\right)dz,
\end{multline*}
where $\Omega=\{(x,y)\in\mathbb{R}^{m+k}: x_{1},y_{1}>1\}$.
\end{cor}
The method can be applied analogously to obtain sharp remainder terms of inequalities \cite[Theorem 3.4, 3.7, 3.8]{yener2018several}.

\subsection{Sharp remainder and stability of the $L^{p}$-Heisenberg-Pauli-Weyl inequality}\label{sec hpw sharp} As noted earlier, the motivation for obtaining the sharp remainder formula of the HPW inequality comes from the question of stability. That is, as shown in \cite{cazacu2024caffarelli}, if we know the explicit form of the \say{Heisenberg deficit}
\begin{align*}
\delta_{2}(f):=\left(\int_{\mathbb{R}^{n}}|\nabla f|^{2}dz\right)\left(\int_{\mathbb{R}^{n}}|z|^{2}|f|^{2}dz\right)-\frac{n^{2}}{4}\left(\int_{\mathbb{R}^{n}}|f|^{2}dz\right)^{2},
\end{align*}
which (for a slightly different deficit) turns out to be
\begin{align*}%\label{lu hpw}
\left(\int_{\mathbb{R}^{n}}|\nabla f|^{2}dz\right)^{\frac{1}{2}}\left(\int_{\mathbb{R}^{n}}|z|^{2}|f|^{2}dz\right)^{\frac{1}{2}}-\frac{n}{2}\int_{\mathbb{R}^{n}}|f|^{2}dz=\frac{\alpha^{2}}{2}\int_{\mathbb{R}^{n}}\left|e^{-\frac{|z|^{2}}{2\alpha^{2}}}\nabla\left(\frac{f}{e^{-\frac{|z|^{2}}{2\alpha^{2}}}}\right)\right|^{2}dz,
\end{align*}
where $\alpha=\left(\frac{\int_{\mathbb{R}^{n}}|f|^{2}|z|^{2}dz}{\int_{\mathbb{R}^{n}}|\nabla f|^{2}dz}\right)^{\frac{1}{4}}$, then with the combination of the optimal Poincar\'e inequality of Gaussian type measure, one can have a certain result relating to the stability of the HPW inequality.

In this context, we apply the main result (Theorem \ref{main thm}) to obtain the sharp remainder term and stability of the $L^{p}$-HPW inequality. In particular, if we set
\begin{align*}
v=\frac{\rho^{\gamma p}}{|x|^{\gamma p}}, \quad 
\phi=e^{-\alpha \rho^{p'}}, \quad 
p'=\frac{p}{p-1} \quad \text{with} \quad \alpha>0,
\end{align*}
we obtain the following identity:
\begin{multline*}%\label{hpw eq discuss}
\int_{\mathbb{R}^n}\frac{\rho^{\gamma p}}{|x|^{\gamma p}}|\nabla_{\gamma}f|^{p}dz = \left(\frac{\alpha p}{p-1}\right)^{p-1}\int_{\mathbb{R}^n}\left(Q-\alpha p\rho^{\frac{p}{p-1}}\right)|f|^{p}dz\\+\int_{\mathbb{R}^n}\frac{\rho^{\gamma p}}{|x|^{\gamma p}}C_{p}\left(\nabla_{\gamma}f,e^{-\alpha \rho^{\frac{p}{p-1}}}\nabla_{\gamma}\left(\frac{f}{e^{-\alpha \rho^{\frac{p}{p-1}}}}\right)\right)dz,
\end{multline*}
which after maximization with respect to $\alpha$ implies the next result.

\begin{cor}\label{cor hpw sharp}
Let $1<p<\infty$ and $p'=\frac{p}{p-1}$. Then, for all non-zero complex-valued $f\in C_{0}^{\infty}(\mathbb{R}^{m+k}\backslash\{x=0\})$, we have
\begin{multline}\label{hpw sharp eq}
\left(\int_{\mathbb{R}^{n}}\frac{\rho^{\gamma p}}{|x|^{\gamma p}}|\nabla_{\gamma}f|^{p}dz\right)\left(\int_{\mathbb{R}^{n}}\rho^{p'}|f|^{p}dz\right)^{p-1}=\left(\frac{Q}{p}\right)^{p}\left(\int_{\mathbb{R}^{n}}|f|^{p}dz\right)^{p}\\+\left(\int_{\mathbb{R}^{n}}\rho^{p'}|f|^{p}dz\right)^{p-1}\int_{\mathbb{R}^{n}}\frac{\rho^{\gamma p}}{|x|^{\gamma p}}C_{p}\left(\nabla_{\gamma}f,e^{-\alpha_{p,Q,\gamma} \rho^{p'}}\nabla_{\gamma}\left(\frac{f}{e^{-\alpha_{p,Q,\gamma} \rho^{p'}}}\right)\right)dz,
\end{multline}
where $\alpha_{p,Q,\gamma}=\frac{Q}{p}\frac{p-1}{p}\frac{\int_{\mathbb{R}^{n}}|f|^{p}dz}{\int_{\mathbb{R}^{n}}\rho^{p'}|f|^{p}dz}$.
\end{cor}

\begin{rem}
We note that (due to Lemmata \ref{lem1}, \ref{lem2} and \ref{lem3}) the identity (\ref{hpw sharp eq}) implies that all extremizers of the $L^{p}$-HPW inequality (related to the Baouendi-Grushin operator) are in the set $E_{p,\gamma}=\{ce^{-\beta\rho^{p'}}:c\in\mathbb{C},\beta>0\}$ and similar remarks can be made for the $\gamma=0$ and $p=2$ cases. 
\end{rem}

When $\gamma=0$, in (\ref{hpw sharp eq}), we obtain the sharp remainder formula of the $L^{p}$-HPW inequality for standard Euclidean gradient.

\begin{cor}\label{cor hpw sharp gamma = 0}
Let $1<p<\infty$ and $p'=\frac{p}{p-1}$. Then, for all non-zero complex-valued $f\in C_{0}^{\infty}(\mathbb{R}^{n})$, we have
\begin{multline}\label{in sec p hpw identity}
\left(\int_{\mathbb{R}^{n}}|\nabla f|^{p}dz\right)\left(\int_{\mathbb{R}^{n}}|z|^{p'}|f|^{p}dz\right)^{p-1}-\left(\frac{n}{p}\right)^{p}\left(\int_{\mathbb{R}^{n}}|f|^{p}dz\right)^{p}
\\=\left(\int_{\mathbb{R}^{n}}|z|^{p'}|f|^{p}dz\right)^{p-1}\int_{\mathbb{R}^{n}}C_{p}\left(\nabla f,e^{-\alpha_{p,n}|z|^{p'}}\nabla\left(\frac{f}{e^{-\alpha_{p,n}|z|^{p'}}}\right)\right)dz,
\end{multline}
where $\alpha_{p,n}=\frac{n}{p}\frac{p-1}{p}\frac{\int_{\mathbb{R}^{n}}|f|^{p}dz}{\int_{\mathbb{R}^{n}}|z|^{p'}|f|^{p}dz}$.
\end{cor}

\begin{rem}
In Corollary \ref{cor hpw sharp} we excluded the singular point $x=0$ from the domain. However, when $\gamma=0$ the singular weight $\tfrac{1}{|x|^{\gamma p}}$ vanishes, and consequently the point $x=0$ can be included in the domain of $f$ without restriction.
\end{rem}

Now let us define the $L^{p}$-Heisenberg deficit:
\begin{align*}
\delta_{p}(f):=\left(\int_{\mathbb{R}^{n}}|\nabla f|^{p}dz\right)\left(\int_{\mathbb{R}^{n}}|z|^{p'}|f|^{p}dz\right)^{p-1}-\left(\frac{n}{p}\right)^{p}\left(\int_{\mathbb{R}^{n}}|f|^{p}dz\right)^{p}.
\end{align*}

Utilizing the sharp remainder formula of the $L^{p}$-HPW inequality (\ref{in sec p hpw identity}) and Lemma \ref{lem1} with the weighted $L^{p}$-Poincar\'e inequality for the log-concave probability measure by Do, Flynn, Lam and Lu \cite[Lemma 4.1]{do2023p}, we are able to obtain the quantitative stability result of the $L^{p}$-HPW inequality.
\begin{thm}\label{stability hpw}
Let $2\leq p<n$ and $p'=\frac{p}{p-1}$. Then, for all non-zero complex-valued $f\in C^{\infty}_{0}(\mathbb{R}^n)$, we have
\begin{align}\label{stab eq}
\delta_{p}(f)\geq \widetilde{C}(n,p)\left(\int_{\mathbb{R}^{n}}|z|^{p'}|f|^{p}dz\right)^{p-1}d(f,E_{p})^{p},
\end{align}
where $\widetilde{C}(n,p)>0$ and $d(f,E_{p})$ is the $L^{p}$ distance to the manifold of extremizers, i.e.
\begin{align*}
d(f,E_{p})^p:=\inf_{c\in\mathbb{C},\beta>0}\left\{\norm{f-ce^{-\beta|z|^{p'}}}^{p}_{L^{p}(\mathbb{R}^n)}\right\}.
\end{align*}
\end{thm}

We note that the $p=2$ case of (\ref{in sec p hpw identity}) gives the sharp remainder formula for the classical $L^{2}$-HPW inequality.

\begin{cor}
For all non-zero complex-valued $f\in C_{0}^{\infty}(\mathbb{R}^{n})$, we have
\begin{multline*}
\left(\int_{\mathbb{R}^{n}}|\nabla f|^{2}dz\right)\left(\int_{\mathbb{R}^{n}}|z|^{2}|f|^{2}dz\right)-\frac{n^2}{4}\left(\int_{\mathbb{R}^{n}}|f|^{2}dz\right)^{2}
\\=\left(\int_{\mathbb{R}^{n}}|z|^{2}|f|^{2}dz\right)\int_{\mathbb{R}^{n}}C_{2}\left(\nabla f,e^{-\alpha_{2,n}|z|^{2}}\nabla\left(\frac{f}{e^{-\alpha_{2,n}|z|^{2}}}\right)\right)dz
\\=\left(\int_{\mathbb{R}^{n}}|z|^{2}|f|^{2}dz\right)\int_{\mathbb{R}^{n}}\left|e^{-\alpha_{2,n}|z|^{2}}\nabla\left(\frac{f}{e^{-\alpha_{2,n}|z|^{2}}}\right)\right|^{2}dz,
\end{multline*}
where $\alpha_{2,n}=\frac{n}{4}\frac{\int_{\mathbb{R}^{n}}|f|^{2}dz}{\int_{\mathbb{R}^{n}}|z|^{2}|f|^{2}dz}$.
\end{cor}

In a similar way, from (\ref{stab eq}), we obtain the stability of the classical $L^{2}$-HPW inequality:
\begin{cor}
For all non-zero complex-valued $f\in C^{\infty}_{0}(\mathbb{R}^n)$, we have
\begin{align*}
\delta_{2}(f)\geq \widetilde{C}(n)\left(\int_{\mathbb{R}^n}|z|^{2}|f|^{2}dz\right)d(f,E)^{2},
\end{align*}
where $\widetilde{C}(n)>0$ and $d(f,E)$ is the $L^{2}$ distance to the manifold of extremizers, i.e
\begin{align*}
d(f,E)^{2}:=\inf_{c\in\mathbb{C},\beta>0}\left\{\norm{f-ce^{\beta|z|^{2}}}^{2}_{L^{2}(\mathbb{R}^n)}\right\}.
\end{align*}
\end{cor}

\subsection{Non-existence of general solutions to parabolic equations}\label{sec: non-exist}

According to Theorem \ref{main thm}, let $v\in C^{1}(\Omega)$ and $w \in L^1_{loc}(\Omega)$ be non-negative weights. If the equation 
\begin{equation}
-\nabla_{\gamma} \cdot \left( v |\nabla_{\gamma} \phi|^{p-2} \nabla_{\gamma} \phi \right) = w \phi^{p-1}
\end{equation}
admits a positive distributional solution, then the following sharp inequality holds:
\begin{equation}
\int_{\Omega} v |\nabla_{\gamma} f|^p \, dz \geq \int_{\Omega} w |f|^p \, dz.
\end{equation}
The sharpness of this inequality implies that if $w$ is replaced by $\tilde{w} = \beta w$ for any $\beta > 1$, the inequality fails. In this section, we discuss how the sharp constants of inequalities involving singular potentials play a crucial role in proving the non-existence of general positive local solutions to nonlinear parabolic equations.

Following the framework established in \cite{Kom06}, we define the notion of a general positive local solution as follows:

\begin{defn}[Positive local solution continuous off of $\mathcal K$]\label{def:2.3}
By a \emph{positive local solution continuous off of $\mathcal K$}, we mean:
\begin{enumerate}
\item[(i)] $\mathcal K$ is a closed Lebesgue null subset of $\Omega$,
\item[(ii)] $u : [0,T) \to L^{1}(\Omega)$ is continuous for some $T>0$,
\item[(iii)] $(z,t) \mapsto u(z,t) \in C\big((\Omega\setminus \mathcal K)\times(0,T)\big)$,
\item[(iv)] $u(z,t)>0$ on $(\Omega\setminus \mathcal K)\times(0,T)$,
\item[(v)] $\lim_{t\to 0} u(\cdot,t)=u_{0}$ in the sense of distributions,
\item[(vi)] $\nabla_{\gamma}u\in L^{2}_{\mathrm{loc}}(\Omega)$
      [resp.\ $\nabla_{\gamma}u\in L^{p}_{\mathrm{loc}}(\Omega)$]
      and $u$ is a solution in the sense of distributions of the PDE.
\end{enumerate}
\end{defn}

\begin{rem}\label{rem:2.4}
If $0<a<b<T$ and $\mathcal K_{0}$ is a compact subset of $\Omega\setminus\mathcal K$,
then
\[
u(z,t)\ge \varepsilon_{1}>0
\quad\text{for } (z,t)\in \mathcal K_{0}\times[a,b]
\]
for some $\varepsilon_{1}>0$.
We can weaken {\rm (iii), (iv)} to be
\begin{enumerate}
\item[(iii)$'$] $u(z,t)$ is positive and locally bounded on $(\Omega\setminus\mathcal K)\times(0,T)$,
\item[(iv)$'$] $\displaystyle \frac{1}{u(z,t)}$ is locally bounded on
$(\Omega\setminus\mathcal K)\times(0,T)$.
\end{enumerate}

If a solution satisfies {\rm (i), (ii), (iii)$'$, (iv)$'$, (v),} and {\rm (vi)},
then we call it a \emph{general positive local solution off of $\mathcal K$}.
This is more general than a positive local solution continuous off of $\mathcal K$.
If $\mathcal K=\emptyset$, we simply call $u$ a
\emph{general positive local solution}.
\end{rem}

\begin{thm}\label{parabolic eq} Let $\frac{2Q}{Q+1} \leq p < 2$, $p-1 < q \leq p-1 + \frac{p}{Q}$, $\lambda\in\mathbb{R}$ and $\Omega \subset \mathbb{R}^N$ denotes a Carnot-Carath\'eodory metric ball. Let $v\in C^{1}(\Omega \backslash \mathcal{K})$ and $\tilde{w}\in L^{1}_{loc}(\Omega\backslash\mathcal K)$ be nonnegative functions, where $\mathcal K$ is a closed Lebesgue null subset of $\Omega$. If
$$
\inf_{0 \neq \phi \in C_c^{\infty}(\Omega \setminus \mathcal K)} \frac{\int_{\Omega} (\epsilon+(1-\epsilon)v(z))|\nabla_\gamma \phi|^p dz - \int_{\Omega} (1-\epsilon) \tilde{w}(z) |\phi|^p dz}{\int_{\Omega} |\phi|^p dz} = -\infty
$$
for some $\epsilon > 0$, then the problem \begin{equation}\label{par. eq}
    \begin{cases}
\frac{\partial u}{\partial t} = \Delta^{v(z)}_{\gamma,p} u + \tilde{w}(z) u^{p-1} + \lambda u^q & \text{in } \Omega \times (0, T), \\
u(z, 0) = u_0(z) \geq 0 & \text{in } \Omega, \\
u(z, t) = 0 & \text{on } \partial \Omega \times (0, T)
\end{cases}
\end{equation} has no generalized positive local solution off of $\mathcal K$. Here, \\$\Delta^{v(z)}_{\gamma,p} u=\nabla_\gamma \cdot\left(v(z)\left|\nabla_\gamma u \right|^{p-2} \nabla_\gamma u\right)$ is the weighted $p$-Grushin operator.
\end{thm}

We now show some applications of Theorem \ref{parabolic eq}, in the next Corollaries. For example, utilizing Corollary \ref{cor log rho}, we can easily get a sharp inequality, for $1<p<\infty$, $\alpha<-1$ and $Q\geq p$, as follows:
$$\int_{B_{R}}^{}\left(\log\frac{R}{\rho}\right)^{\alpha+p}|\nabla_{\gamma}f|^{p}dz\geq\left( \frac{|\alpha+1|}{p} \right) ^{p}\int_{B_{R}}^{}\left( \log \frac{R}{\rho} \right)^{\alpha}\frac{|x|^{\gamma p}}{\rho^{\gamma p+p}}|f|^{p}dz,$$
with \say{virtual} extremizer $\tilde{\phi}=\left(\log \frac{R}{\rho}\right)^{\frac{|\alpha+1|}{p}}$. 

\begin{cor}\label{cor1}
Let $\alpha<-1$, $Q\geq p$ and $p-1 < q \leq p-1 + \frac{p}{Q}$. If $\frac{2Q}{Q+1} \leq p < 2$ and $c>\left( \frac{|\alpha+1|}{p} \right) ^{p}$, then the problem \begin{equation}\label{par. eq with log}
\begin{small}
    \begin{cases}
\frac{\partial u}{\partial t} = \nabla_\gamma \cdot\left(\left(\log\frac{R}{\rho}\right)^{\alpha+p}\left|\nabla_\gamma u\right|^{p-2} \nabla_\gamma u\right) + c\left( \log \frac{R}{\rho} \right)^{\alpha} \frac{|x|^{\gamma p}}{\rho^{\gamma p+p}}u^{p-1} + \lambda u^q & \text{in } {B_{R}} \times (0, T), \\
u(z, 0) = u_0(z) \geq 0 & \text{in } {B_{R}}, \\
u(z, t) = 0 & \text{on } \partial {B_{R}} \times (0, T)
\end{cases}
\end{small}
\end{equation} has no generalized positive local solution off of $\mathcal \{(0,0)\}$.   
\end{cor}

From Corollary \ref{cor nhc}, we have a sharp inequality, for $1<p<\infty$ and $R>0$, as follows:
\begin{align}\label{pde app 2 ineq}
\int_{B_{R}}^{}|\nabla_{\gamma}f|^{p}dz\geq\left( \frac{p-1}{p} \right) ^{p}\int_{B_{R}}^{}\frac{|x|^{\gamma p}}{(R-\rho)^{p}\rho^{\gamma p}}|f|^{p}dz,    
\end{align}
with \say{virtual} extremizer $\tilde{\phi}=(R-\rho)^{\frac{p-1}{p}}$.
The next Corollary shows how the sharp constant in the inequality (\ref{pde app 2 ineq}) plays a role in the proof of the non-existence of a general positive local solution.

\begin{cor}\label{cor2}
Let $1<p<\infty$, $R>0$ and $p-1 < q \leq p-1 + \frac{p}{Q}$. If $\frac{2Q}{Q+1} \leq p < 2$ and $c>\left( \frac{p-1}{p} \right) ^{p}$, then the problem \begin{equation}\label{par. eq with R-r}
\begin{cases}
\frac{\partial u}{\partial t} = \nabla_\gamma \cdot\left(\left|\nabla_\gamma u\right|^{p-2} \nabla_\gamma u\right) + c\frac{|x|^{\gamma p}}{(R-\rho)^{p}\rho^{\gamma p}} u^{p-1} + \lambda u^q & \text{in } {B_{R}} \times (0, T), \\
u(z, 0) = u_0(z) \geq 0 & \text{in } {B_{R}}, \\
u(z, t) = 0 & \text{on } \partial {B_{R}} \times (0, T)
\end{cases}
\end{equation} has no generalized positive local solution.    
\end{cor}

\section{Sharp Remainder of the Hardy-Poincar\'e Inequality With General Weights--Proof of Theorem \ref{main thm}}\label{proof of main}

\begin{proof}[\textbf{Proof of Theorem \ref{main thm}}]
For any complex-valued $f \in C_0^{\infty}(\Omega)$, we set 
\[
\psi := \frac{f}{\phi},
\]
where $\phi$ is a positive solution to
\begin{align}\label{p-grushin proof}
-\nabla_{\gamma} \cdot \left( v \left|\nabla_{\gamma} \phi\right|^{p-2} \nabla_{\gamma} \phi \right) \geq w \phi^{p-1}
\end{align}
in the distributional sense in $\Omega$. Applying the formula for the $C_p$-functional with the vectors 
\[
\xi = \psi \nabla_{\gamma}\phi + \phi\nabla_{\gamma}\psi, \quad 
\eta = \phi \nabla_{\gamma} \psi,
\]
we obtain
\begin{align}\label{step1}
C_{p}(\psi \nabla_{\gamma}\phi + \phi\nabla_{\gamma}\psi, \phi\nabla_{\gamma}\psi) \nonumber
&= \left| \psi \nabla_\gamma \phi + \phi \nabla_\gamma \psi \right|^p 
- \left| \psi \nabla_\gamma \phi \right|^p 
\\& \quad- p \left| \psi \nabla_\gamma \phi \right|^{p-2} 
\text{Re}\!\left( (\psi \nabla_\gamma \phi) \cdot \overline{(\phi \nabla_\gamma \psi)} \right) \nonumber
\\
&= |\nabla_{\gamma}f|^{p} 
- |\nabla_{\gamma}\phi|^{p} |\psi|^{p} \nonumber
\\& \quad- p \text{Re} |\psi|^{p-2} \psi \overline{\nabla_{\gamma}\psi} 
|\nabla_{\gamma}\phi|^{p-2}\phi\nabla_{\gamma}\phi \nonumber
\\
&= |\nabla_{\gamma}f|^{p}
- |\nabla_{\gamma}\phi|^{p} |\psi|^{p} \nonumber
\\& \quad- \nabla_{\gamma}(|\psi|^{p}) 
|\nabla_{\gamma}\phi|^{p-2} \phi \nabla_{\gamma}\phi.
\end{align}
Integrating both sides of (\ref{step1}) and multiplying by $v$, we have
\begin{align*}
\int_{\Omega} v |\nabla_{\gamma}f|^{p} dz
&= \int_{\Omega} 
C_{p}(\psi \nabla_{\gamma}\phi + \phi\nabla_{\gamma}\psi, \phi\nabla_{\gamma}\psi) v dz
+ \int_{\Omega} |\nabla_{\gamma}\phi|^{p} |\psi|^{p} v dz
\\
&\quad + \int_{\Omega} 
\nabla_{\gamma}(|\psi|^{p}) 
|\nabla_{\gamma}\phi|^{p-2} 
\phi \nabla_{\gamma}\phi v dz.
\end{align*}
By integrating by parts, we get
\begin{align*}
\int_{\Omega} v |\nabla_{\gamma}f|^{p} dz
&= \int_{\Omega} 
C_{p}(\psi \nabla_{\gamma}\phi + \phi\nabla_{\gamma}\psi, \phi\nabla_{\gamma}\psi) v dz
+ \int_{\Omega} |\nabla_{\gamma}\phi|^{p} |\psi|^{p} v dz
\\
&\quad - \int_{\Omega} 
|\psi|^{p} \nabla_{\gamma} \cdot 
\left( \phi v |\nabla_{\gamma}\phi|^{p-2} \nabla_{\gamma}\phi \right) dz.
\end{align*}
Expanding the divergence term yields
\begin{align*}
\int_{\Omega} v |\nabla_{\gamma}f|^{p} dz
&= \int_{\Omega} 
C_{p}(\psi \nabla_{\gamma}\phi + \phi\nabla_{\gamma}\psi, \phi\nabla_{\gamma}\psi) v dz
+ \int_{\Omega} |\nabla_{\gamma}\phi|^{p} |\psi|^{p} v dz
\\
&\quad - \int_{\Omega} 
|\psi|^{p} \nabla_{\gamma}\phi 
\left( v |\nabla_{\gamma}\phi|^{p-2}\nabla_{\gamma}\phi \right) dz
 - \int_{\Omega} 
|\psi|^{p} \phi 
\nabla_{\gamma} \cdot \left( v |\nabla_{\gamma}\phi|^{p-2} \nabla_{\gamma}\phi \right) dz.
\end{align*}
Simplifying, we get
\begin{align*}
\int_{\Omega} v |\nabla_{\gamma}f|^{p} dz
&= \int_{\Omega} 
C_{p}(\psi \nabla_{\gamma}\phi + \phi\nabla_{\gamma}\psi, \phi\nabla_{\gamma}\psi) v dz
+ \int_{\Omega} |\nabla_{\gamma}\phi|^{p} |\psi|^{p} v dz
\\
&\quad - \int_{\Omega} 
|\nabla_{\gamma}\phi|^{p} |\psi|^{p} v dz
- \int_{\Omega} 
|\psi|^{p} \phi 
\nabla_{\gamma} \cdot \left( v |\nabla_{\gamma}\phi|^{p-2} \nabla_{\gamma}\phi \right) dz
\\
&= \int_{\Omega} 
C_{p}(\psi \nabla_{\gamma}\phi + \phi\nabla_{\gamma}\psi, \phi\nabla_{\gamma}\psi) v dz
 - \int_{\Omega} 
|\psi|^{p} \phi 
\nabla_{\gamma} \cdot \left( v |\nabla_{\gamma}\phi|^{p-2} \nabla_{\gamma}\phi \right) dz
\\&=\int_{\Omega} 
C_{p}(\psi \nabla_{\gamma}\phi + \phi\nabla_{\gamma}\psi, \phi\nabla_{\gamma}\psi) v dz
+ \int_{\Omega}v|\nabla_{\gamma}\phi|^{p-2}\nabla_{\gamma}\phi\cdot\nabla_{\gamma}\left(|\psi|^{p}\phi\right)dz.
\end{align*}
Since (\ref{p-grushin proof}), by definition, we have
\begin{align}\label{distributional def}
\int_{\Omega}v|\nabla_{\gamma}\phi|^{p-2}\nabla_{\gamma}\phi\cdot\nabla_{\gamma}\varphi dz \geq \int_{\Omega}w\phi^{p-1}\varphi dz
\end{align}
for all $\varphi \in C_{0}^{\infty}(\Omega)$. Applying (\ref{distributional def}) for $\varphi=|\psi|^{p}\phi$, we obtain
\begin{align*}
\int_{\Omega} v |\nabla_{\gamma}f|^{p} dz
&\geq \int_{\Omega} 
C_{p}(\psi \nabla_{\gamma}\phi + \phi\nabla_{\gamma}\psi, \phi\nabla_{\gamma}\psi) v dz
+ \int_{\Omega} |\psi|^{p} w |\phi|^{p} dz.
\end{align*}
After performing the change of variables $\psi=\frac{f}{\phi}$, we get
\begin{align*}
\int_{\Omega} v |\nabla_{\gamma}f|^{p} dz
\geq \int_{\Omega} 
C_{p}\left( \nabla_{\gamma}f, \phi \nabla_{\gamma} \left( \frac{f}{\phi} \right) \right) v dz
+ \int_{\Omega} w|f|^{p} dz.
\end{align*}
The proof is now complete.
\end{proof}

\section{Proofs of Corollary \ref{cor rec KY} and of Applications of the Main Result}\label{proofs of cor}

\begin{proof}[\textbf{Proof of Corollary \ref{cor rec KY}}]
Using Lemma \ref{lem1}, in (\ref{main result}), we immediately get (\ref{kom1}) with the explicit constant $c_1(p)$:
\begin{align*}
\int_{\Omega} v\left|\nabla_\gamma f\right|^p dz \geq \int_{\Omega} w|f|^p dz+c_1(p) \int_{\Omega} v\left|\nabla_\gamma \left(\frac{f}{\phi}\right)\right|^p \phi^p dz.
\end{align*}
To get (\ref{kom2}), we set
\begin{align}\label{variables}
\xi=\nabla_{\gamma}f, \quad \eta=\phi\nabla_{\gamma}\left(\frac{f}{\phi}\right) \quad \text{and} \quad \xi-\eta=\frac{f}{\phi}\nabla_{\gamma}\phi
\end{align}
so that
\begin{align*}
|\xi|+|\xi-\eta|=|\eta+(\xi-\eta)|+|\xi-\eta|\leq 2|\xi-\eta|+|\eta|\leq 2\left(|\xi-\eta|+|\eta|\right).
\end{align*}
Thus, we have
\begin{align*}
\left(|\xi| + |\xi - \eta|\right)^{2-p}
\leq
\left[2\left(|\xi - \eta| + |\eta|\right)\right]^{2-p}
=
2^{2-p}\left(|\xi - \eta| + |\eta|\right)^{2-p},
\end{align*}
which gives us
\begin{align}\label{xi eta est}
\frac{|\eta|^2}{\left(|\xi| + |\xi - \eta|\right)^{2-p}}\geq
\frac{|\eta|^2}{2^{2-p}\left(|\xi - \eta| + |\eta|\right)^{2-p}}
=
2^{p-2}
\frac{|\eta|^2}{\left(\left|\frac{f}{\phi}\nabla\phi\right| + \left|\phi\nabla \left(\frac{f}{\phi}\right)\right|\right)^{2-p}}.
\end{align}
Substituting (\ref{variables}) back to (\ref{xi eta est}), we derive
\begin{align}\label{xi eta last}
\frac{\left| \phi \nabla_{\gamma} \left( \frac{f}{\phi} \right) \right|^2}{\left( |\nabla_{\gamma}f| + \left| \nabla_{\gamma}f - \phi \nabla_{\gamma} \left( \frac{f}{\phi} \right) \right| \right)^{2-p}} \geq 2^{p-2}\frac{\left|\phi\nabla_{\gamma}\left(\frac{f}{\phi}\right)\right|^{2}}{\left(\left|\frac{f}{\phi}\nabla\phi\right| + \left|\phi\nabla \left(\frac{f}{\phi}\right)\right|\right)^{2-p}}.
\end{align}
Using Lemma \ref{lem2}, in (\ref{main result}), we get
\begin{multline*}
\int_{\Omega}^{}v|\nabla_{\gamma}f|^{p}dz \geq \int_{\Omega}^{}w|f|^{p}dz + \int_{\Omega}^{}vC_{p}\left( \nabla_{\gamma}f,\phi\nabla_{\gamma}\left( \frac{f}{\phi} \right)  \right) dz
\\ \geq \int_{\Omega}^{}w|f|^{p}dz+c_{2}(p)\int_{\Omega}^{}\frac{v\left| \phi \nabla_{\gamma} \left( \frac{f}{\phi} \right) \right|^2}{\left( |\nabla_{\gamma}f| + \left| \nabla_{\gamma}f - \phi \nabla_{\gamma} \left( \frac{f}{\phi} \right) \right| \right)^{2-p}}dz.
\end{multline*}
Now we utilize the obtained estimate (\ref{xi eta last}):
\begin{align*}
\int_{\Omega}^{}v|\nabla_{\gamma}f|^{p}dz\geq \int_{\Omega}^{}w|f|^{p}dz+\widetilde{c_{2}}(p)\int_{\Omega}^{}\frac{v\left| \phi \nabla_{\gamma} \left( \frac{f}{\phi} \right) \right|^2}{\left(\left|\frac{f}{\phi}\nabla\phi\right| + \left|\phi\nabla \left(\frac{f}{\phi}\right)\right|\right)^{2-p}}dz,
\end{align*}
where $\widetilde{c_{2}}(p)=2^{p-2}c_2(p)$.
\end{proof}

\subsection{Proofs of sharp remainder formulas of weighted Hardy-type and \\Poincar\'e inequalities--Proofs of Corollaries \ref{cor nhc}-\ref{cor several yener}}

\begin{proof}[\textbf{Proof of Corollary \ref{cor nhc}}]
Let
\begin{align*}
v=1, \quad \phi=(R-\rho)^{\frac{p-1}{p}}.
\end{align*}
We need to calculate 
\begin{align}\label{w1}
w=\frac{-\nabla_{\gamma}\cdot \left( v|\nabla_{\gamma}\phi|^{p-2}\nabla_{\gamma}\phi \right) }{\phi^{p-1}}=\frac{-\nabla_{\gamma}\cdot T}{\phi^{p-1}}=\frac{\Phi}{\phi^{p-1}},
\end{align}
where
\begin{align}\label{Phi1}
T=v|\nabla_{\gamma}\phi|^{p-2}\nabla_{\gamma}\phi \quad \text{ and } \quad \Phi=-\nabla_{\gamma}\cdot T.
\end{align}
Now since $\phi=\phi(\rho)$, we have
\begin{align*}
\nabla_{\gamma}\phi=\phi'(\rho)\nabla_{\gamma}\rho \quad \text{and} \quad |\nabla_{\gamma}\phi|=|\phi'(\rho)||\nabla_{\gamma}\rho|,
\end{align*}
which gives us
\begin{align*}
\nabla_{\gamma}\phi=-c(R-\rho)^{-c-1}\nabla_{\gamma}\rho \quad \text{and} \quad |\nabla_{\gamma}\phi|=c\frac{|\nabla_{\gamma}\rho|}{\left( R-\rho \right)^{\frac{1}{p}} }, 
\end{align*}
where $c=\frac{p-1}{p}$. Therefore, we have
\begin{align}\label{T1}
T=-c^{p-1}\frac{|\nabla_{\gamma}\rho|^{p-2}\nabla_{\gamma}\rho}{\left( R-\rho \right) ^{c}}.   
\end{align}
Substituting (\ref{T1}) to (\ref{Phi1}), we get
\begin{align*}
\Phi&=c^{p-1}\nabla_{\gamma}\cdot \left( \frac{|\nabla_{\gamma}\rho|^{p-2}\nabla_{\gamma}\rho}{\left( R-\rho \right) ^{c}} \right) 
\\&= c^{p-1}\left[ \nabla_{\gamma}\left( \left( R-\rho \right) ^{-c} \right)|\nabla_{\gamma}\rho|^{p-2}\nabla_{\gamma}\rho + \left( R-\rho \right) ^{-c}\Delta_{\gamma}\rho  \right] 
\\&=c^{p-1}\left[ c\left( R-\rho \right) ^{-c-1}|\nabla_{\gamma}\rho|^{p}+\frac{Q-1}{\rho}\frac{|\nabla_{\gamma}\rho|^{p}}{(R-\rho)^{c}} \right] 
\\&=c^{p}\frac{|\nabla_{\gamma}\rho|^{p}}{(R-\rho)^{c+1}}+(Q-1)c^{p-1}\frac{|\nabla_{\gamma}\rho|^{p}}{\rho(R-\rho)^{c}}.
\end{align*}
Putting (\ref{Phi1}) to (\ref{w1}), we derive
\begin{align*}
w&=c^{p}\frac{|\nabla_{\gamma}\rho|^{p}}{(R-\rho)^{pc+1}}+\frac{|\nabla_{\gamma}\rho|^{p}}{\rho(R-\rho)^{pc}}
\\&=\left( \frac{p-1}{p} \right) ^{p}\frac{|x|^{\gamma p}}{(R-\rho)^{p}\rho^{\gamma p}}+(Q-1)\left( \frac{p-1}{p} \right) ^{p-1}\frac{|x|^{\gamma p}}{(R-\rho)^{p-1}\rho^{\gamma p + 1}}
\end{align*}
and the proof is complete.
\end{proof}

\begin{proof}[\textbf{Proof of Corollary \ref{cor dam}}]
Let
\begin{align*}
v=|x|^{\beta-\gamma p}\rho^{(1+\gamma)p-\alpha}, \qquad \phi=\rho^{-\frac{Q+\beta-\alpha}{p}}.
\end{align*}
We need to calculate
\begin{align}\label{w2}
w=\frac{-\nabla_{\gamma}\cdot\left(v|\nabla_{\gamma}\phi|^{p-2}\nabla_{\gamma}\phi\right)}{\phi^{p-1}}
=\frac{-\nabla_{\gamma}\cdot T}{\phi^{p-1}}
=\frac{\Phi}{\phi^{p-1}},
\end{align}
where
\begin{align}\label{Phi2}
T=v|\nabla_{\gamma}\phi|^{p-2}\nabla_{\gamma}\phi \quad \text{and} \quad \Phi=-\nabla_{\gamma}\cdot T.
\end{align}
Since $\phi=\phi(\rho)$, we have
\begin{align*}
\nabla_{\gamma}\phi=\phi'(\rho)\nabla_{\gamma}\rho \quad \text{and} \quad |\nabla_{\gamma}\phi|=|\phi'(\rho)||\nabla_{\gamma}\rho|.
\end{align*}
Let $c=\dfrac{Q+\beta-\alpha}{p}$, then, this yields
\begin{align*}
\nabla_{\gamma}\phi=-c\rho^{-c-1}\nabla_{\gamma}\rho \quad \text{and} \quad
|\nabla_{\gamma}\phi|=c\frac{|\nabla_{\gamma}\rho|}{\rho^{c+1}}.
\end{align*}
Therefore,
\begin{align}
T&=-c^{p-1}|x|^{\beta-\gamma p}\rho^{(1+\gamma)p-\alpha}
\frac{|\nabla_{\gamma}\rho|^{p-2}}{\rho^{(c+1)(p-2)}}
\rho^{-c-1}\nabla_{\gamma}\rho\notag
\\&=-c^{p-1}\frac{|x|^{\beta-\gamma p}}{\rho^{(c+1)(p-1)+\alpha-(1+\gamma)p}}|\nabla_{\gamma}\rho|^{p-2}\nabla_{\gamma}\rho\notag
\\&=-c^{p-1}\frac{|x|^{\beta-\gamma p}}{\rho^{(c+1)(p-1)+\alpha-(1+\gamma)p}}
\frac{|x|^{\gamma(p-2)}}{\rho^{\gamma(p-2)}}\nabla_{\gamma}\rho\notag
\\&=-c^{p-1}\frac{|x|^{\beta-2\gamma}}{\rho^{(c+1)(p-1)+\alpha-p-2\gamma}}\nabla_{\gamma}\rho.\label{T2}
\end{align}
Substituting \eqref{T2} into \eqref{Phi2}, we get
\begin{align*}
\Phi
&=c^{p-1}\nabla_{\gamma}\cdot\left(|x|^{\beta-2\gamma}\rho^{-(c+1)(p-1)-\alpha+p+2\gamma}\nabla_{\gamma}\rho\right)
\\&=c^{p-1}\left(Q+\beta-2\gamma-(c+1)(p-1)-\alpha+p+2\gamma-1\right)\frac{|x|^{\beta}}{\rho^{1+(c+1)(p-1)+\alpha-p}}
\\&=c^{p-1}\left(Q+\beta-\alpha-cp+c\right)
\frac{|x|^{\beta}}{\rho^{c(p-1)+\alpha}}
\\&=c^{p}\frac{|x|^{\beta}}{\rho^{c(p-1)+\alpha}}.
\end{align*}
Putting \eqref{Phi2} into \eqref{w2}, we obtain
\begin{align*}
w
&=c^{p}\frac{|x|^{\beta}}{\rho^{\alpha}}
=\left(\frac{Q+\beta-\alpha}{p}\right)^{p}\frac{|x|^{\beta}}{\rho^{\alpha}},
\end{align*}
completing the proof.
\end{proof}

\begin{proof}[\textbf{Proof of Corollary \ref{cor dam |x|}}]
Let
\begin{align*}
v=|x|^{\alpha+p}, \qquad \phi=|x|^{s}, \qquad s=\frac{|m+\alpha|}{p}.
\end{align*}
We need to calculate
\begin{align}\label{w3}
w=\frac{-\nabla_{\gamma}\cdot\left(v|\nabla_{\gamma}\phi|^{p-2}\nabla_{\gamma}\phi\right)}{\phi^{p-1}}
=\frac{-\nabla_{\gamma}\cdot T}{\phi^{p-1}}
=\frac{\Phi}{\phi^{p-1}},
\end{align}
where
\begin{align*}
T=v|\nabla_{\gamma}\phi|^{p-2}\nabla_{\gamma}\phi, \qquad \Phi=-\nabla_{\gamma}\cdot T.
\end{align*}
Since $\phi=\phi(x)$ depends only on the $x$ variables, we have
\begin{align*}
\nabla_{\gamma}\phi=(\nabla_{x}\phi,0), \quad \nabla_{x}\phi=s|x|^{s-2}x, \quad |\nabla_{\gamma}\phi|=s|x|^{s-1}.
\end{align*}
Therefore,
\begin{align*}
|\nabla_{\gamma}\phi|^{p-2}=s^{p-2}|x|^{(s-1)(p-2)}.
\end{align*}
Hence
\begin{align*}
T&=|x|^{\alpha+p}s^{p-2}|x|^{(s-1)(p-2)}\left(s|x|^{s-2}x,0\right)
\\&=s^{p-1}|x|^{\alpha+p+(s-1)(p-2)+(s-2)}(x,0).
\end{align*}
A simple computation gives
\begin{align*}
\alpha+p+(s-1)(p-2)+(s-2)=\alpha+s(p-1),
\end{align*}
so
\begin{align*}
T=s^{p-1}|x|^{\alpha+s(p-1)}(x,0).
\end{align*}
Taking the divergence, we get
\begin{align}\label{Phi3}
\Phi=-\sum_{i=1}^{m}\partial_{x_{i}}\left(s^{p-1}|x|^{\alpha+s(p-1)}x_{i}\right)
&=-s^{p-1}\left(m+\alpha+s(p-1)\right)|x|^{\alpha+s(p-1)}.
\end{align}
Since $\phi^{p-1}=|x|^{s(p-1)}$, from \eqref{w3} and \eqref{Phi3}, we obtain
\begin{align*}
w&=-s^{p-1}\left(m+\alpha+s(p-1)\right)\frac{|x|^{\alpha+s(p-1)}}{|x|^{s(p-1)}}.
\end{align*}
Using $-(m+\alpha)=|m+\alpha|$, we deduce
\begin{align*}
w=s^{p}|x|^{\alpha}=\left(\frac{|m+\alpha|}{p}\right)^{p}|x|^{\alpha}.
\end{align*}
The proof is complete.
\end{proof}

\begin{proof}[\textbf{Proof of Corollary \ref{cor poin}}]
Take
\begin{align}\label{pair}
v=1 \quad \text{and} \quad \phi=\phi_1:=\text{the first eigenfunction of } -\Delta_{\gamma,p}.
\end{align}
Assuming that $\phi_1>0$, (\ref{pair}) immediately gives us
\begin{align*}
w=\frac{-\nabla_{\gamma}\cdot\left(|\nabla_{\gamma}\phi_{1}|^{p-2}\nabla_{\gamma}\phi_{1}\right)}{\phi_{1}^{p-1}}=\frac{- \Delta_{\gamma,p} \phi_{1}}{\phi^{p-1}_{1}}=\lambda_{1},
\end{align*}
where $\lambda_{1}$ is a number (eigenvalue) such that
\begin{align*}
- \Delta_{\gamma,p} \phi_{1} = \lambda_{1} \phi_{1}^{p-1}.
\end{align*}
The proof is complete.
\end{proof}

\begin{proof}[\textbf{Proof of Corollary \ref{cor log rho}}]
Let
\begin{align*}
v=\left( \log \frac{R}{\rho} \right)^{\alpha+p}, \quad 
\phi=\left( \log \frac{R}{\rho} \right)^{s}, \quad 
s=\frac{|\alpha+1|}{p}.
\end{align*}
We need to calculate
\begin{align}\label{wlog1}
w=\frac{-\nabla_{\gamma}\cdot\left(v|\nabla_{\gamma}\phi|^{p-2}\nabla_{\gamma}\phi\right)}{\phi^{p-1}}
=\frac{-\nabla_{\gamma}\cdot T}{\phi^{p-1}}
=\frac{\Phi}{\phi^{p-1}},
\end{align}
where
\begin{align}\label{Philog1}
T=v|\nabla_{\gamma}\phi|^{p-2}\nabla_{\gamma}\phi, \quad 
\Phi=-\nabla_{\gamma}\cdot T.
\end{align}
Since $\phi=\phi(\rho)$, we have
\begin{align*}
\nabla_{\gamma}\phi=\phi'(\rho)\nabla_{\gamma}\rho, \quad 
|\nabla_{\gamma}\phi|=|\phi'(\rho)||\nabla_{\gamma}\rho|.
\end{align*}
Now
\begin{align*}
\phi'(\rho) 
= \frac{d}{d\rho}\left( \log \frac{R}{\rho} \right)^{s} 
= -s\frac{\left(\log \frac{R}{\rho}\right)^{s-1}}{\rho},
\end{align*}
so that
\begin{align*}
\nabla_{\gamma}\phi 
=-s\frac{\left(\log \frac{R}{\rho}\right)^{s-1}}{\rho}\nabla_{\gamma}\rho, \quad 
|\nabla_{\gamma}\phi|=s\frac{\left(\log \frac{R}{\rho}\right)^{s-1}}{\rho}|\nabla_{\gamma}\rho|.
\end{align*}
Therefore
\begin{align}
T&=-\left(\log \frac{R}{\rho}\right)^{\alpha+p}
\left(s\frac{\left(\log \frac{R}{\rho}\right)^{s-1}}{\rho}|\nabla_{\gamma}\rho|\right)^{p-2}
s\frac{\left(\log \frac{R}{\rho}\right)^{s-1}}{\rho}\nabla_{\gamma}\rho\notag
\\&=-s^{p-1}\left(\log \frac{R}{\rho}\right)^{\alpha+p+(s-1)(p-1)}
\frac{|\nabla_{\gamma}\rho|^{p-2}\nabla_{\gamma}\rho}{\rho^{p-1}}\notag
\\&=-s^{p-1}L^{\alpha+1+s(p-1)}\rho^{-(p-1)}|\nabla_{\gamma}\rho|^{p-2}\nabla_{\gamma}\rho,\label{Tlog1}
\end{align}
where
\begin{align*}
L=\log \frac{R}{\rho}.
\end{align*}
From \eqref{Philog1} and \eqref{Tlog1}, we get
\begin{align*}
\Phi&=-\nabla_{\gamma}\cdot T
=s^{p-1}\nabla_{\gamma}\cdot\left(L^{\alpha+1+s(p-1)}\rho^{-(p-1)}|\nabla_{\gamma}\rho|^{p-2}\nabla_{\gamma}\rho\right)
\\&=s^{p-1}\nabla_{\gamma}A(\rho)|\nabla_{\gamma}\rho|^{p-2}\nabla_{\gamma}\rho
+s^{p-1}A(\rho)\frac{Q-1}{\rho}|\nabla_{\gamma}\rho|^{p},
\end{align*}
where
\begin{align*}
A(\rho)=L^{\alpha+1+s(p-1)}\rho^{-(p-1)}.
\end{align*}
We compute
\begin{align*}
\nabla_{\gamma}A(\rho)
&=\nabla_{\gamma}\left(L^{\alpha+1+s(p-1)}\right)\rho^{-(p-1)}
+L^{\alpha+1+s(p-1)}\nabla_{\gamma}\rho^{-(p-1)}
\\&=-(\alpha+1+s(p-1))L^{\alpha+s(p-1)}\rho^{-p}\nabla_{\gamma}\rho
-(p-1)L^{\alpha+1+s(p-1)}\rho^{-p}\nabla_{\gamma}\rho
\\&=-L^{\alpha+s(p-1)}\rho^{-p}\nabla_{\gamma}\rho\left[(\alpha+1+s(p-1))+(p-1)L\right].
\end{align*}
Since $\alpha+1=-|\alpha+1|$, we obtain
\begin{align*}
\nabla_{\gamma}A(\rho)
=L^{\alpha+s(p-1)}\rho^{-p}\nabla_{\gamma}\rho\left[s-(p-1)L\right].
\end{align*}
Thus, we get
\begin{align}
\Phi&=s^{p-1}L^{\alpha+s(p-1)}\rho^{-p}\left[s-(p-1)L\right]|\nabla_{\gamma}\rho|^{p}
+s^{p-1}(Q-1)L^{\alpha+1+s(p-1)}\frac{|\nabla_{\gamma}\rho|^{p}}{\rho^{p}}\notag
\\&=s^{p}L^{\alpha+s(p-1)}\frac{|\nabla_{\gamma}\rho|^{p}}{\rho^{p}}
+s^{p-1}(Q-p)L^{\alpha+1+s(p-1)}\frac{|\nabla_{\gamma}\rho|^{p}}{\rho^{p}}.\label{Philog2}
\end{align}
Dividing \eqref{Philog2} by $\phi^{p-1}=L^{s(p-1)}$ and using \eqref{wlog1}, we obtain
\begin{align*}
w=s^{p}L^{\alpha}\frac{|\nabla_{\gamma}\rho|^{p}}{\rho^{p}}
+s^{p-1}(Q-p)L^{\alpha+1}\frac{|\nabla_{\gamma}\rho|^{p}}{\rho^{p}}.
\end{align*}
Putting everything back, we get exactly 
\begin{align*}
w=\left(\frac{|\alpha+1|}{p}\right)^p\left(\log\frac{R}{\rho}\right)^{\alpha}\frac{|x|^{\gamma p}}{\rho^{\gamma p + p}}+\left(\frac{|\alpha+p|}{p}\right)^{p-1}(Q-p)\left(\log\frac{R}{\rho}\right)^{\alpha+1}\frac{|x|^{\gamma p}}{\rho^{\gamma p +p}}
\end{align*}
and the proof is complete.
\end{proof}

\begin{proof}[\textbf{Proof of Corollary \ref{cor log |x|}}]
Let
\begin{align*}
v=\left(\log\frac{R}{|x|}\right)^{\alpha+p}, \quad 
\phi=\left(\log\frac{R}{|x|}\right)^{s}, \quad
s=\frac{|\alpha+1|}{p},
\end{align*}
and set
\begin{align*}
L=\log\frac{R}{|x|}.
\end{align*}
We need to calculate
\begin{align}\label{wxlog1}
w=\frac{-\nabla_{\gamma}\cdot\left(v|\nabla_{\gamma}\phi|^{p-2}\nabla_{\gamma}\phi\right)}{\phi^{p-1}}
=\frac{-\nabla_{\gamma}\cdot T}{\phi^{p-1}}
=\frac{\Phi}{\phi^{p-1}},
\end{align}
where
\begin{align}\label{Phixlog1}
T=v|\nabla_{\gamma}\phi|^{p-2}\nabla_{\gamma}\phi, \quad 
\Phi=-\nabla_{\gamma}\cdot T.
\end{align}
Since $\phi=\phi(|x|)$ depends only on the $x$ variables, we have
\begin{align*}
\nabla_{\gamma}\phi=(\nabla_x\phi,0), \quad 
\nabla_x\phi=sL^{s-1}\nabla_x\left(\log\frac{R}{|x|}\right)
=-sL^{s-1}\frac{x}{|x|^{2}}.
\end{align*}
Thus, we have
\begin{align*}
|\nabla_{\gamma}\phi|=s\frac{L^{s-1}}{|x|}, \quad 
|\nabla_{\gamma}\phi|^{p-2}=s^{p-2}L^{(s-1)(p-2)}|x|^{-(p-2)}.
\end{align*}
Therefore, we get
\begin{align}
T&=v|\nabla_{\gamma}\phi|^{p-2}\nabla_{\gamma}\phi\notag
\\&=-s^{p-1}L^{\alpha+p+(s-1)(p-2)+(s-1)}|x|^{-(p-2)}\frac{x}{|x|^{2}}\notag
\\&=-s^{p-1}L^{\alpha+1+s(p-1)}|x|^{-p}(x,0).\label{Txlog1}
\end{align}
Let
\begin{align*}
A(|x|)=L^{\alpha+1+s(p-1)}|x|^{-p}.
\end{align*}
Since $\nabla_x\cdot(Ax)=mA+|x|A'(|x|)$, from \eqref{Phixlog1} and \eqref{Txlog1} we obtain
\begin{align*}
\Phi=-\nabla_{\gamma}\cdot T=s^{p-1}\left(mA+|x|A'(|x|)\right).
\end{align*}
A direct computation gives
\begin{align*}
|x|A'(|x|)=|x|^{-p}\left(-pL^{k}-kL^{k-1}\right), \quad 
k=\alpha+1+s(p-1),
\end{align*}
and hence
\begin{align}
\Phi=s^{p-1}|x|^{-p}\left((m-p)L^{k}-kL^{k-1}\right).\label{Phixlog2}
\end{align}
Dividing \eqref{Phixlog2} by $\phi^{p-1}=L^{s(p-1)}$ and using \eqref{wxlog1}, we get
\begin{align*}
w=s^{p-1}|x|^{-p}\left((m-p)L^{\alpha+1}-(\alpha+1+s(p-1))L^{\alpha}\right).
\end{align*}
Since $\alpha+1=-|\alpha+1|$ we have $s=\frac{|\alpha+1|}{p}$, which yields
\begin{align*}
w&=\left(\frac{|\alpha+1|}{p}\right)^{p}\frac{L^{\alpha}}{|x|^{p}}
+\left(\frac{|\alpha+1|}{p}\right)^{p-1}(m-p)\frac{L^{\alpha+1}}{|x|^{p}}
\\&=\left(\frac{|\alpha+1|}{p}\right)^{p}\frac{\left(\log \frac{R}{|x|}\right)^{\alpha}}{|x|^{p}}+\left(\frac{|\alpha+1|}{p}\right)^{p-1}(m-p)\frac{\left(\log \frac{R}{|x|}\right)^{\alpha+1}}{|x|^{p}}
\end{align*}
and the proof is complete.
\end{proof}

\begin{proof}[\textbf{Proof of Corollary \ref{cor rho alpha}}]
Let
\begin{align*}
v=\left( 1+\rho^{a} \right)^{\alpha(p-1)}, \quad 
\phi=\left( 1+\rho^{a} \right)^{1-\alpha}, \quad 
a=\frac{p}{p-1}, \quad \alpha>1.
\end{align*}
Set $b=\alpha-1>0$, then
\begin{align*}
\phi=\left( 1+\rho^{a} \right)^{-b}.
\end{align*}
We need to calculate
\begin{align}\label{w new}
w=\frac{-\nabla_{\gamma}\cdot \left( v|\nabla_{\gamma}\phi|^{p-2}\nabla_{\gamma}\phi \right)}{\phi^{p-1}}
=\frac{-\nabla_{\gamma}\cdot T}{\phi^{p-1}}
=\frac{\Phi}{\phi^{p-1}},
\end{align}
where
\begin{align}\label{Phi new}
T=v|\nabla_{\gamma}\phi|^{p-2}\nabla_{\gamma}\phi, \quad \Phi=-\nabla_{\gamma}\cdot T.
\end{align}
Since $\phi=\phi(\rho)$, we have
\begin{align*}
\phi'(\rho)=-b\,a\,\rho^{a-1}(1+\rho^{a})^{-b-1}, \quad
\nabla_{\gamma}\phi=\phi'(\rho)\nabla_{\gamma}\rho, \quad
|\nabla_{\gamma}\phi|=|\phi'(\rho)||\nabla_{\gamma}\rho|.
\end{align*}
Hence
\begin{align*}
|\phi'|^{p-2}\phi'=-(ba)^{p-1}\rho^{(a-1)(p-1)}(1+\rho^{a})^{-(b+1)(p-1)}.
\end{align*}
Multiplying by $v=(1+\rho^{a})^{(b+1)(p-1)}$ gives
\begin{align*}
T&=-(ba)^{p-1}\rho^{(a-1)(p-1)}|\nabla_{\gamma}\rho|^{p-2}\nabla_{\gamma}\rho.
\end{align*}
Since $(a-1)(p-1)=1$, we obtain
\begin{align}\label{T new}
T=-(ba)^{p-1}\rho\,|\nabla_{\gamma}\rho|^{p-2}\nabla_{\gamma}\rho.
\end{align}
Using the radial identity
\begin{align*}
\nabla_{\gamma}\cdot\left(A(\rho)|\nabla_{\gamma}\rho|^{p-2}\nabla_{\gamma}\rho\right)
=\left(A'(\rho)+\frac{Q-1}{\rho}A(\rho)\right)|\nabla_{\gamma}\rho|^{p}
\end{align*}
with $A(\rho)=\rho$ and $A'=1$, we get from \eqref{Phi new} and \eqref{T new}
\begin{align*}
\Phi=(ba)^{p-1}Q|\nabla_{\gamma}\rho|^{p}.
\end{align*}
Dividing by $\phi^{p-1}=(1+\rho^{a})^{-b(p-1)}$ and using \eqref{w new}, we have
\begin{align*}
w=(ba)^{p-1}Q(1+\rho^{a})^{b(p-1)}|\nabla_{\gamma}\rho|^{p}.
\end{align*}
Recalling that $|\nabla_{\gamma}\rho|^{p}=\dfrac{|x|^{\gamma p}}{\rho^{\gamma p}}$ and $ba=\dfrac{(\alpha-1)p}{p-1}$, we arrive at
\begin{align*}
w=Q\left(\frac{(\alpha-1)p}{p-1}\right)^{p-1}
\left(1+\rho^{\frac{p}{p-1}}\right)^{(\alpha-1)(p-1)}
\frac{|x|^{\gamma p}}{\rho^{\gamma p}}.
\end{align*}
The proof is complete.
\end{proof}

\begin{proof}[\textbf{Proof of Corollary \ref{cor super}}]
Let
\begin{align*}
v=\frac{(a+b\rho^{\alpha})^{\beta}}{\rho^{\ell p}}, \quad 
\phi=\rho^{-s}, \quad s=\frac{Q-p\ell-p}{p}.
\end{align*}
We need to calculate
\begin{align}\label{wrad1}
w=\frac{-\nabla_{\gamma}\cdot\left(v|\nabla_{\gamma}\phi|^{p-2}\nabla_{\gamma}\phi\right)}{\phi^{p-1}}
=\frac{-\nabla_{\gamma}\cdot T}{\phi^{p-1}}
=\frac{\Phi}{\phi^{p-1}},
\end{align}
where
\begin{align}\label{Phirad1}
T=v|\nabla_{\gamma}\phi|^{p-2}\nabla_{\gamma}\phi, \quad \Phi=-\nabla_{\gamma}\cdot T.
\end{align}
Since $\phi=\phi(\rho)$, we have
\begin{align*}
\nabla_{\gamma}\phi=\phi'(\rho)\nabla_{\gamma}\rho=-s\rho^{-s-1}\nabla_{\gamma}\rho, \quad
|\nabla_{\gamma}\phi|=s\rho^{-s-1}|\nabla_{\gamma}\rho|.
\end{align*}
Therefore
\begin{align}
T=-s^{p-1}(a+b\rho^{\alpha})^{\beta}\rho^{-\ell p-(s+1)(p-1)}|\nabla_{\gamma}\rho|^{p-2}\nabla_{\gamma}\rho.\label{Trad1}
\end{align}
Let
\begin{align*}
A(\rho)=(a+b\rho^{\alpha})^{\beta}\rho^{-\ell p-(s+1)(p-1)}.
\end{align*}
Using the radial identity
\begin{align*}
\nabla_{\gamma}\cdot\left(A(\rho)|\nabla_{\gamma}\rho|^{p-2}\nabla_{\gamma}\rho\right)
=\left(A'(\rho)+\frac{Q-1}{\rho}A(\rho)\right)|\nabla_{\gamma}\rho|^{p},
\end{align*}
we obtain from \eqref{Phirad1} and \eqref{Trad1}
\begin{align*}
\Phi=s^{p-1}\left(A'(\rho)+\frac{Q-1}{\rho}A(\rho)\right)|\nabla_{\gamma}\rho|^{p}.
\end{align*}
A direct computation gives
\begin{align*}
\frac{A'(\rho)}{A(\rho)}
=\frac{\beta\alpha b\,\rho^{\alpha-1}}{a+b\rho^{\alpha}}
-\frac{\ell p+(s+1)(p-1)}{\rho},
\end{align*}
and thus
\begin{align*}
A'(\rho)+\frac{Q-1}{\rho}A(\rho)
=A(\rho)\left[\frac{\beta\alpha b\,\rho^{\alpha-1}}{a+b\rho^{\alpha}}
+\frac{Q-1-\ell p-(s+1)(p-1)}{\rho}\right].
\end{align*}
Since $s=\frac{Q-p\ell-p}{p}$, we have
\begin{align*}
Q-1-\ell p-(s+1)(p-1)=s.
\end{align*}
Therefore
\begin{align}
\Phi=s^{p-1}A(\rho)|\nabla_{\gamma}\rho|^{p}
\left[\frac{\beta\alpha b\,\rho^{\alpha-1}}{a+b\rho^{\alpha}}+\frac{s}{\rho}\right].\label{Phirad2}
\end{align}
Dividing \eqref{Phirad2} by $\phi^{p-1}=\rho^{-s(p-1)}$ and using \eqref{wrad1}, we get
\begin{align*}
w=s^{p-1}A(\rho)\rho^{s(p-1)}|\nabla_{\gamma}\rho|^{p}
\left[\frac{\beta\alpha b\,\rho^{\alpha-1}}{a+b\rho^{\alpha}}+\frac{s}{\rho}\right].
\end{align*}
Since $A(\rho)\rho^{s(p-1)}=(a+b\rho^{\alpha})^{\beta}\rho^{-\ell p-(p-1)}$ and $|\nabla_{\gamma}\rho|^{p}=\frac{|x|^{\gamma p}}{\rho^{\gamma p}}$, it follows that
\begin{align*}
w&=s^{p}(a+b\rho^{\alpha})^{\beta}\frac{|x|^{\gamma p}}{\rho^{\gamma p+\ell p+p}}
+s^{p-1}\beta\alpha b(a+b\rho^{\alpha})^{\beta-1}\frac{|x|^{\gamma p}}{\rho^{\gamma p+\ell p+p-\alpha}}
\\&=\left(\frac{Q-p\ell-p}{p}\right)^{p}(a+b\rho^{\alpha})^{\beta}\frac{|x|^{\gamma p}}{\rho^{\gamma p+\ell p+p}}
\\&\quad+\left(\frac{Q-p\ell-p}{p}\right)^{p-1}\beta\alpha b(a+b\rho^{\alpha})^{\beta-1}\frac{|x|^{\gamma p}}{\rho^{\gamma p+\ell p+p-\alpha}}.
\end{align*}
The proof is complete.    
\end{proof}

\begin{proof}[\textbf{Proof of Corollary \ref{cor several yener}}]
Let
\begin{align*}
v=\left(\frac{y_{1}}{|x|^{\gamma}}\right)^{p-2}\log x_{1}, \qquad \phi=\log y_{1},
\end{align*}
with $x_{1},y_{1}>1$. We need to calculate
\begin{align}\label{w newlogy}
w=\frac{-\nabla_{\gamma}\cdot\left(v|\nabla_{\gamma}\phi|^{p-2}\nabla_{\gamma}\phi\right)}{\phi^{p-1}}
=\frac{-\nabla_{\gamma}\cdot T}{\phi^{p-1}}
=\frac{\Phi}{\phi^{p-1}},
\end{align}
where
\begin{align}\label{Phi newlogy}
T=v|\nabla_{\gamma}\phi|^{p-2}\nabla_{\gamma}\phi, \qquad \Phi=-\nabla_{\gamma}\cdot T.
\end{align}
We have
\begin{align*}
X_{i}(\log y_{1})=\frac{\partial}{\partial x_{i}}\log y_{1}=0, \quad i=1,\dots,n,
\end{align*}
and
\begin{align*}
Y_{j}(\log y_{1})=\left(|x|^{\gamma}\frac{\partial}{\partial y_{j}}\right)(\log y_{1})
=\begin{cases}
\dfrac{|x|^{\gamma}}{y_{1}}, & j=1,\\
0, & j\neq 1.
\end{cases}
\end{align*}
Hence
\begin{align*}
\nabla_{\gamma}\log y_{1}=\left(0,\dots,0,\frac{|x|^{\gamma}}{y_{1}},0,\dots,0\right),
\end{align*}
and therefore
\begin{align*}
|\nabla_{\gamma}\log y_{1}|^{p-2}=\left(\frac{|x|^{\gamma}}{y_{1}}\right)^{p-2}.
\end{align*}
It follows from \eqref{Phi newlogy} that
\begin{align*}
T&=\left(\frac{y_{1}}{|x|^{\gamma}}\right)^{p-2}\log x_{1}\left(\frac{|x|^{\gamma}}{y_{1}}\right)^{p-2}\left(0,\dots,0,\frac{|x|^{\gamma}}{y_{1}},0,\dots,0\right)\\
&=\left(0,\dots,0,\frac{|x|^{\gamma}\log x_{1}}{y_{1}},0,\dots,0\right).
\end{align*}
Taking the divergence, we obtain
\begin{align*}
\Phi&=-Y_{1}\left(\frac{|x|^{\gamma}\log x_{1}}{y_{1}}\right)
=-|x|^{2\gamma}\log x_{1}\frac{\partial}{\partial y_{1}}\left(\frac{1}{y_{1}}\right)
=\frac{|x|^{2\gamma}\log x_{1}}{y_{1}^{2}}.
\end{align*}
Using \eqref{w newlogy} and $\phi^{p-1}=\log^{p-1}y_{1}$, we get
\begin{align*}
w=\frac{|x|^{2\gamma}\log x_{1}}{y_{1}^{2}\log^{p-1}y_{1}},
\end{align*}
which completes the proof.
\end{proof}

\subsection{Proof of the sharp remainder formula of the Heisenberg-Pauli-Weyl inequality--Proof of Corollary \ref{cor hpw sharp}}

\begin{proof}[\textbf{Proof of Corollary \ref{cor hpw sharp}}]
Let
\begin{align*}
v=\frac{\rho^{\gamma p}}{|x|^{\gamma p}}, \quad 
\phi=e^{-\alpha \rho^{p'}}, \quad 
p'=\frac{p}{p-1}, \quad \alpha>0.
\end{align*}
We need to calculate
\begin{align}\label{w hpw}
w=\frac{-\nabla_{\gamma}\cdot\left(v|\nabla_{\gamma}\phi|^{p-2}\nabla_{\gamma}\phi\right)}{\phi^{p-1}}
=\frac{-\nabla_{\gamma}\cdot T}{\phi^{p-1}}
=\frac{\Phi}{\phi^{p-1}},
\end{align}
where
\begin{align}\label{Phi hpw}
T=v|\nabla_{\gamma}\phi|^{p-2}\nabla_{\gamma}\phi, \qquad \Phi=-\nabla_{\gamma}\cdot T.
\end{align}
Since $\phi=\phi(\rho)$, we have
\begin{align*}
\phi'(\rho)=-\alpha p'\,\rho^{p'-1}e^{-\alpha\rho^{p'}}, \quad
\nabla_{\gamma}\phi=\phi'(\rho)\nabla_{\gamma}\rho, \quad
|\nabla_{\gamma}\phi|=|\phi'(\rho)||\nabla_{\gamma}\rho|.
\end{align*}
Thus
\begin{align*}
|\nabla_{\gamma}\phi|^{p-2}\nabla_{\gamma}\phi
=-(\alpha p')^{p-1}\rho^{(p'-1)(p-1)}e^{-\alpha(p-1)\rho^{p'}}
|\nabla_{\gamma}\rho|^{p-2}\nabla_{\gamma}\rho,
\end{align*}
and therefore
\begin{align}\label{T hpw}
T=-(\alpha p')^{p-1}A(\rho)v|\nabla_{\gamma}\rho|^{p-2}\nabla_{\gamma}\rho,
\end{align}
where
\begin{align*}
A(\rho)=\rho^{(p'-1)(p-1)}e^{-\alpha(p-1)\rho^{p'}}.
\end{align*}
Using the product rule together with the radial identity
\begin{align*}
\nabla_{\gamma}\cdot\big(|\nabla_{\gamma}\rho|^{p-2}\nabla_{\gamma}\rho\big)
=\frac{Q-1}{\rho}|\nabla_{\gamma}\rho|^{p},
\end{align*}
we get from \eqref{T hpw} and \eqref{Phi hpw}
\begin{align*}
\Phi&=(\alpha p')^{p-1}\nabla_{\gamma}\cdot\Big[A(\rho)v|\nabla_{\gamma}\rho|^{p-2}\nabla_{\gamma}\rho\Big]
\\&=(\alpha p')^{p-1}\Big[\nabla_{\gamma}\left(A(\rho) v\right)|\nabla_{\gamma}\rho|^{p-2}\nabla_{\gamma}\rho+A(\rho)v\Delta_{\gamma,p}\rho\Big]
\\&=(\alpha p')^{p-1}\Big[A'(\rho)v|\nabla_{\gamma}\rho|^{p}+A(\rho)|\nabla_{\gamma}\rho|^{p-2}\nabla_{\gamma}v\cdot\nabla_{\gamma}\rho+A(\rho)v\Delta_{\gamma,p}\rho\Big]
\\&=(\alpha p')^{p-1}\left[ A(\rho)|\nabla_{\gamma}\rho|^{p-2}\nabla_{\gamma}v\cdot\nabla_{\gamma}\rho
+v\left(A'(\rho)+\tfrac{Q-1}{\rho}A(\rho)\right)|\nabla_{\gamma}\rho|^{p}\right]
\\&=(\alpha p')^{p-1}v\left(A'(\rho)+\tfrac{Q-1}{\rho}A(\rho)\right)|\nabla_{\gamma}\rho|^{p}.
\end{align*}
Therefore, we have
\begin{align}\label{phi hpw 2}
\Phi=(\alpha p')^{p-1}\left(A'(\rho)+\frac{Q-1}{\rho}A(\rho)\right).
\end{align}
Since
\begin{align*}
A'(\rho)=\rho^{(p'-1)(p-1)-1}e^{-\alpha(p-1)\rho^{p'}}
\Big((p'-1)(p-1)-\alpha(p-1)p'\rho^{p'}\Big),
\end{align*}
we have
\begin{multline*}
A'(\rho)+\frac{Q-1}{\rho}A(\rho)
\\=\rho^{(p'-1)(p-1)-1}e^{-\alpha(p-1)\rho^{p'}}
\Big((p'-1)(p-1)+Q-1-\alpha(p-1)p'\rho^{p'}\Big).
\end{multline*}
Dividing by $\phi^{p-1}=e^{-\alpha(p-1)\rho^{p'}}$ and using (\ref{phi hpw 2}) with \eqref{w hpw}, we get
\begin{align*}
w=(\alpha p')^{p-1}\rho^{(p'-1)(p-1)-1}
\Big((p'-1)(p-1)+Q-1-\alpha(p-1)p'\rho^{p'}\Big).
\end{align*}
With $p'=\frac{p}{p-1}$ we have $(p'-1)(p-1)=1$, hence
\begin{align*}
w=(\alpha p')^{p-1}\Big(Q-\alpha(p-1)p'\rho^{p'}\Big)
=\left(\frac{\alpha p}{p-1}\right)^{p-1}\left(Q-\alpha p\rho^{\frac{p}{p-1}}\right).
\end{align*}
By Theorem \ref{main thm}, this gives us
\begin{multline}\label{hpw eq}
\int_{\mathbb{R}^n}\frac{\rho^{\gamma p}}{|x|^{\gamma p}}|\nabla_{\gamma}f|^{p}dz = \left(\frac{\alpha p}{p-1}\right)^{p-1}\int_{\mathbb{R}^n}\left(Q-\alpha p\rho^{\frac{p}{p-1}}\right)|f|^{p}dz\\+\int_{\mathbb{R}^n}\frac{\rho^{\gamma p}}{|x|^{\gamma p}}C_{p}\left(\nabla_{\gamma}f,e^{-\alpha \rho^{\frac{p}{p-1}}}\nabla_{\gamma}\left(\frac{f}{e^{-\alpha \rho^{\frac{p}{p-1}}}}\right)\right)dz.
\end{multline}
Now set
\begin{align}\label{hpw change}
K:=\int_{\mathbb{R}^n}\frac{\rho^{\gamma p}}{|x|^{\gamma p}}|\nabla_{\gamma}f|^{p}dz, \quad L:=\int_{\mathbb{R}^n}\rho^{p'}|f|^{p}dz, \quad M:=\int_{\mathbb{R}^n}|f|^{p}dz.
\end{align}
Taking (\ref{hpw change}) and (\ref{hpw eq}) into consideration, we get 
\begin{align}\label{k=}
K=\left(\frac{\alpha p}{p-1}\right)^{p-1}(QM-\alpha pL)+R(\alpha),
\end{align}
where
\begin{align*}
R(\alpha)=\int_{\mathbb{R}^n}\frac{\rho^{\gamma p}}{|x|^{\gamma p}}C_{p}\left(\nabla_{\gamma}f,e^{-\alpha \rho^{p'}}\nabla_{\gamma}\left(\frac{f}{e^{-\alpha \rho^{p'}}}\right)\right)dz.
\end{align*}
Now let
\begin{align}\label{s=alp p}
s:=\alpha p
\end{align}
and with (\ref{s=alp p}) set
\begin{align}\label{Gs}
G(s):=\left(\frac{s}{p-1}\right)^{p-1}(QM-sL)=(p-1)^{-(p-1)}(s^{p-1}QM-s^{p}L).
\end{align}
Next, we maximize (\ref{Gs}) with respect to $s$:
\begin{multline*}
G'(s)=(p-1)^{-(p-1)}s^{p-2}((p-1)QM-psL)=0 \\ \iff s=0 \quad \text{or} \quad s=\frac{(p-1)QM}{pL}.
\end{multline*}
Since $1<p<\infty$ and $\alpha>0$, the only possible critical point is
\begin{align}\label{s crit}
s=\frac{(p-1)QM}{pL} \iff \alpha=\frac{(p-1)QM}{p^{2}L}:=\alpha_{p,Q,\gamma}
\end{align}
for all non-zero complex-valued $f\in C^{\infty}_{0}(\mathbb{R}^n)$. Substituting (\ref{s crit}) to (\ref{Gs}), we get
\begin{align*}
G(s)=\left(\frac{QM}{pL}\right)^{p-1}\frac{QM}{p}=\left(\frac{Q}{p}\right)^{p}\frac{M^{p}}{L^{p-1}}.
\end{align*}
Therefore, from (\ref{k=}), we have
\begin{align}\label{pre last hpw}
K=\left(\frac{Q}{p}\right)^{p}\frac{M^{p}}{L^{p-1}}+R(\alpha_{p,Q,\gamma}).
\end{align}
Multiplying both sides of (\ref{pre last hpw}) by $L^{p-1}$, we derive
\begin{align*}
KL^{p-1}=\left(\frac{Q}{p}\right)^{p}M^{p}+L^{p-1}R(\alpha_{p,Q,\gamma}).
\end{align*}
Thus, from (\ref{hpw change}), we conclude
\begin{multline*}
\left(\int_{\mathbb{R}^n}\frac{\rho^{\gamma p}}{|x|^{\gamma p}}|\nabla_{\gamma}f|^{p}dz\right)\left(\int_{\mathbb{R}^n}\rho^{p'}|f|^{p}dz\right)^{p-1}=\left(\frac{Q}{p}\right)^{p}\left(\int_{\mathbb{R}^n}|f|^{p}dz\right)^{p}\\+\left(\int_{\mathbb{R}^n}\rho^{p'}|f|^{p}dz\right)^{p-1}\int_{\mathbb{R}^n}\frac{\rho^{\gamma p}}{|x|^{\gamma p}}C_{p}\left(\nabla_{\gamma}f,e^{-\alpha_{p,Q,\gamma} \rho^{p'}}\nabla_{\gamma}\left(\frac{f}{e^{-\alpha_{p,Q,\gamma} \rho^{p'}}}\right)\right)dz
\end{multline*}
and the proof is complete.
\end{proof}

\subsection{Proof of the stability of the Heisenberg-Pauli-Weyl inequality--Proof of Theorem \ref{stability hpw}}\label{sec stab hpw}
\begin{proof}[\textbf{Proof of Theorem \ref{stability hpw}}]
Let $0\neq f\in C^{\infty}_{0}(\mathbb{R}^n)$. From Corollary \ref{cor hpw sharp gamma = 0}, we have
\begin{multline}\label{step1 stab}
\left(\int_{\mathbb{R}^{n}}|\nabla f|^{p}dz\right)\left(\int_{\mathbb{R}^{n}}|z|^{p'}|f|^{p}dz\right)^{p-1}-\left(\frac{n}{p}\right)^{p}\left(\int_{\mathbb{R}^{n}}|f|^{p}dz\right)^{p}
\\=\left(\int_{\mathbb{R}^{n}}|z|^{p'}|f|^{p}dz\right)^{p-1}\int_{\mathbb{R}^{n}}C_{p}\left(\nabla f,e^{-\alpha_{p,n}|z|^{p'}}\nabla\left(\frac{f}{e^{-\alpha_{p,n}|z|^{p'}}}\right)\right)dz,
\end{multline}
where $1<p<\infty$ and $p'=\frac{p}{p-1}>0$ and $\alpha_{p,n}=\frac{n}{p}\frac{p-1}{p}\frac{\int_{\mathbb{R}^{n}}|f|^{p}dz}{\int_{\mathbb{R}^{n}}|z|^{p'}|f|^{p}dz}$. Using Lemma $\ref{lem1}$ to (\ref{step1 stab}), we get
\begin{multline}\label{stab main part}
\left(\int_{\mathbb{R}^{n}}|\nabla f|^{p}dz\right)\left(\int_{\mathbb{R}^{n}}|z|^{p'}|f|^{p}dz\right)^{p-1}-\left(\frac{n}{p}\right)^{p}\left(\int_{\mathbb{R}^{n}}|f|^{p}dz\right)^{p}\\\geq c_{1}(p)\left(\int_{\mathbb{R}^{n}}|z|^{p'}|f|^{p}dz\right)^{p-1}\int_{\mathbb{R}^n}\left|e^{-\alpha_{p,n}|z|^{p'}}\nabla\left(\frac{f}{e^{-\alpha_{p,n}|z|^{p'}}}\right)\right|^{p}dz
\\=c_{1}(p)\left(\int_{\mathbb{R}^{n}}|z|^{p'}|f|^{p}dz\right)^{p-1}\int_{\mathbb{R}^n}\left|\nabla\left(\frac{f}{e^{-\alpha_{p,n}|z|^{p'}}}\right)\right|^{p}e^{-p\alpha_{p,n}|z|^{p'}}dz
\end{multline}
for $2 \leq p < \infty$. Let us define
\begin{align*}
u:=\frac{f}{e^{-\alpha_{p,n}|z|^{p'}}}.
\end{align*}
Then, we get that
\begin{align}\label{almost stability}
\int_{\mathbb{R}^n}\left|\nabla\left(\frac{f}{e^{-\alpha_{p,n}|z|^{p'}}}\right)\right|^{p}e^{-p\alpha_{p,n}|z|^{p'}}dz=\int_{\mathbb{R}^n}^{}\left|\nabla u\right|^{p}e^{-p\alpha_{p,n}|z|^{p'}}dz,
\end{align}
where $2\leq p<\infty$, $\alpha_{p,n}>0$ and $p'=\frac{p}{p-1}>1$. Now let us recall the weighted $L^{p}$-Poincar\'e inequality for the log-concave probability measure by Do, Flynn, Lam and Lu \cite[Lemma 4.1]{do2023p}: for some $\delta>0$, $n-p>\mu\geq0$ and $\alpha\geq\frac{n-p-\mu}{n-p}$, we have for $u\in C^{\infty}_{0}(\mathbb{R}^n\backslash\{0\})$ that
\begin{align}\label{lu poincare}
\int_{\mathbb{R}^n}\frac{|\nabla u|^{p}}{|z|^{\mu}}e^{-\delta|z|^{\alpha}}dz\geq C(n,p,\alpha,\delta,\mu)\inf_{c}\int_{\mathbb{R}^n}\frac{|u-c|^{p}}{|z|^{\frac{n\mu}{n-p}}}e^{-\delta|z|^{\alpha}}dz,
\end{align}
where $C(N,p,\alpha,\delta,\mu)>0$ is some universal positive constant. Setting $\mu=0$, $\delta=p \alpha_{p,n}>0$, $\alpha=p'>1$, in (\ref{lu poincare}), we get that for $u\in C^{\infty}_{0}(\mathbb{R}^n)$
\begin{align*}
\int_{\mathbb{R}^n}|\nabla u|^{p}e^{-p \alpha_{p,n}|z|^{p'}}dz\geq C(n,p)\inf_{c}\int_{\mathbb{R}^n}|u-c|^{p}e^{-p \alpha_{p,n}|z|^{p'}}dz,
\end{align*}
which when applied to (\ref{almost stability}) gives
\begin{align}\label{stab last}
\int_{\mathbb{R}^n}\left|\nabla\left(\frac{f}{e^{-\alpha_{p,n}|z|^{p'}}}\right)\right|^{p}e^{-p\alpha_{p,n}|z|^{p'}}dz \nonumber
&=\int_{\mathbb{R}^n}^{}\left|\nabla u\right|^{p}e^{-p\alpha_{p,n}|z|^{p'}}dz \nonumber
\\&\geq C(n,p)\inf_{c}\int_{\mathbb{R}^n}|u-c|^{p}e^{-p \alpha_{p,n}|z|^{p'}}dz \nonumber
\\&=C(n,p)\inf_{c}\int_{\mathbb{R}^n}\left|\frac{f}{e^{-\alpha_{p,n}|z|^{p'}}}-c\right|^{p}e^{-p \alpha_{p,n}|z|^{p'}}dz \nonumber
\\&=C(n,p)\inf_{c}\int_{\mathbb{R}^n}\left|f-ce^{- \alpha_{p,n}|z|^{p'}}\right|^{p}dz.
\end{align}
Using (\ref{stab last}) in (\ref{stab main part}), we finally obtain
\begin{multline*}
\left(\int_{\mathbb{R}^{n}}|\nabla f|^{p}dz\right)\left(\int_{\mathbb{R}^{n}}|z|^{p'}|f|^{p}dz\right)^{p-1}-\left(\frac{n}{p}\right)^{p}\left(\int_{\mathbb{R}^{n}}|f|^{p}dz\right)^{p}\\\geq c_{1}(p)C(n,p)\left(\int_{\mathbb{R}^{n}}|z|^{p'}|f|^{p}dz\right)^{p-1}\inf_{c}\int_{\mathbb{R}^n}\left|f-ce^{- \alpha_{p,n}|z|^{p'}}\right|^{p}dz
\\=\widetilde{C}(n,p)\left(\int_{\mathbb{R}^{n}}|z|^{p'}|f|^{p}dz\right)^{p-1}\inf_{c}\int_{\mathbb{R}^n}\left|f-ce^{- \alpha_{p,n}|z|^{p'}}\right|^{p}dz.
\end{multline*}
\end{proof}

\subsection{Proofs of the non-existence of general positive local solutions--Proofs of Theorem \ref{parabolic eq} and Corollaries \ref{cor1}-\ref{cor2}}

The proof of Theorem \ref{parabolic eq} utilizes \cite[Proposition 2.2]{Kom06}, which we recall below:
\begin{prop}\label{prop:2.2}
Let $\Omega = B(z_0,r)$ be a $d_C$-metric ball in $\mathbb{R}^{m+k}$, let
$1 < p < Q$, let $M \in L^{Q/p}(\Omega)$, and let
$\phi \in W^{1,p}_0(\Omega)$ have compact support in $\Omega$.
Then for each $\varepsilon > 0$ there exists a positive constant
$C(\varepsilon)$ such that
\begin{equation}\label{eq:2.7}
\int_\Omega M(z)\,|\phi|^p\,dz
\leq \frac{\varepsilon}{1-\varepsilon}
\int_\Omega |\nabla_\gamma \phi|^p\,dz
+ C(\varepsilon)\int_\Omega |\phi|^p\,dz.
\end{equation}
\end{prop}

\begin{proof}[Proof of Theorem \ref{parabolic eq}]
We proceed by contradiction. Suppose that for any $T > 0$, there exists a generalized positive local solution $u : [0, T) \to L^1(\Omega)$ to \eqref{par. eq} in the domain $(\Omega \setminus \mathcal{K}) \times (0, T)$, where $u_0 \geq 0$ and $u_0 \not\equiv 0$. Multiplying the first equation of \eqref{par. eq} with $|\phi|^p/u^{p-1}$ for $\phi \in C_{0}^{\infty}(\Omega \setminus \mathcal{K})$ and integrating over $\Omega$ yields
\begin{multline}\label{mult. phi}
\frac{1}{2-p} \frac{d}{d t} \int_{\Omega} u^{2-p} |\phi|^p dz = \int_{\Omega} \Delta^{v(z)}_{\gamma,p} u \left( \frac{|\phi|^p}{u^{p-1}} \right) dz + \int_{\Omega} \tilde{w}(z) |\phi|^p dz\\+ \int_{\Omega} \lambda u^{q-p+1} |\phi|^p dz.
\end{multline}
Then, the divergence theorem gives
\begin{equation}\label{int by parts 1}
\int_{\Omega} \Delta^{v(z)}_{\gamma,p} u \left(\frac{|\phi|^p}{u^{p-1}}\right) dz
= - \int_{\Omega} v(z)|\nabla_\gamma u|^{p-2} \nabla_\gamma u \cdot \nabla_\gamma \left(\frac{|\phi|^p}{u^{p-1}}\right) dz.
\end{equation}
A direct computation shows that
\begin{equation}\label{int by parts 2}
v(z)|\nabla_\gamma u|^{p-2}\nabla_\gamma u \cdot \nabla_\gamma\!\left(\frac{|\phi|^p}{u^{p-1}}\right)
= p v(z)|\nabla_\gamma u|^{p-2} \frac{|\phi|^{p-1}}{u^{p-1}} \nabla_\gamma u \cdot \nabla_\gamma|\phi|
- (p-1)v(z)\frac{|\phi|^p}{u^p} |\nabla_\gamma u|^p.
\end{equation}
Now we substitute \eqref{int by parts 2} into \eqref{int by parts 1} to get
\begin{multline}\label{int p lap}
\int_{\Omega} \Delta^{v(z)}_{\gamma,p} u \left(\frac{|\phi|^p}{u^{p-1}}\right) dz
= (p-1) \int_{\Omega} v(z)|\nabla_\gamma u|^p \frac{|\phi|^p}{u^p} dz
\\- p \int_{\Omega} v(z)|\nabla_\gamma u|^{p-2} \frac{|\phi|^{p-1}}{u^{p-1}} \nabla_\gamma u \cdot \nabla_\gamma|\phi| dz.
\end{multline}
Applying the Cauchy-Schwarz inequality to the second integral on the right-hand side of \eqref{int p lap}, we obtain
\begin{multline}\label{int p lap 2}
\int_{\Omega} \Delta^{v(z)}_{\gamma,p} u \left(\frac{|\phi|^p}{u^{p-1}}\right) dz
\geq (p-1) \int_{\Omega} v(z)|\nabla_\gamma u|^p \frac{|\phi|^p}{u^p} dz
\\- p \int_{\Omega} v(z)|\nabla_\gamma u|^{p-1}|\nabla_\gamma \phi|\frac{|\phi|^{p-1}}{u^{p-1}} dz.
\end{multline}
We now invoke the following elementary inequality: for $p>1$ and distinct positive real numbers $a$ and $b$, 
\[
a^p - b^p - p b^{p-1}(a - b) > 0,
\]
which implies that
\[
(p-1) b^p - p b^{p-1}a > - a^p.
\]
We can take $b = \left|v(z)^\frac{1}{p}\frac{\phi}{u}\nabla_\gamma u\right|$ and $a = |v(z)^\frac{1}{p}\nabla_\gamma \phi|$; then we have
\[
(p-1) \int_{\Omega} v(z)|\nabla_\gamma u|^p \frac{|\phi|^p}{u^p} dz
- p \int_{\Omega} v(z)|\nabla_\gamma u|^{p-1}|\nabla_\gamma \phi| \frac{|\phi|^{p-1}}{u^{p-1}} dz
\geq - \int_{\Omega} v(z)|\nabla_\gamma \phi|^p dz.
\]
Hence, we conclude using in \eqref{int p lap 2} that
\begin{equation} \label{int p lap 3}
\int_{\Omega} \Delta^{v(z)}_{\gamma,p} u \left(\frac{|\phi|^p}{u^{p-1}}\right) dz
\geq - \int_{\Omega} v(z)|\nabla_\gamma \phi|^p dz,
\end{equation}
provided $\phi \neq 0$. Substituting \eqref{int p lap 3} into \eqref{mult. phi} and integrating from $t_1$ to $t_2$ , where $0 < t_1 < t_2 < T$ , we obtain

\begin{equation}\label{<L}
\int_{\Omega} \tilde{w}(z) |\phi|^p dz - \int_{\Omega} v(z)|\nabla_\gamma \phi|^p dz \leq L,
\end{equation}
where
\begin{multline}\label{L}
L := \frac{1}{(2-p)(t_2 - t_1)} \int_{\Omega} (u^{2-p}(z, t_2) - u^{2-p}(z, t_1)) |\phi|^p dz
\\- \lambda \int_{t_1}^{t_2} \int_{\Omega} u^{q-p+1} |\phi|^p dz dt.
\end{multline}
Using Jensen's inequality for concave functions, we obtain, since $\frac{2Q}{Q+1} \leq p < 2$,
$$
\int_{\Omega} (u^{2-p}(z, t_i))^{\frac{Q}{p}} dz \leq C(|\Omega|) \left( \int_{\Omega} u(z, t_i) dz \right)^{\frac{(2-p)Q}{p}} < \infty.
$$
Therefore
$$
u^{2-p}(z, t_i) \in L^{Q/p}(\Omega).
$$
Next, we examine the second integral on the right-hand side of \eqref{L}. By defining $F(z) := \lambda \int_{t_1}^{t_2} u(z,t)^{q-p+1} dt$ and applying Jensen's inequality, it follows that
$$
F(z) \in L^{Q/p}(\Omega).
$$
As a direct consequence of Proposition \ref{prop:2.2}, we obtain
\begin{equation}\label{L<}
L \leq \frac{\epsilon}{1-\epsilon} \int_{\Omega} |\nabla_\gamma \phi|^p dz + C(\epsilon) \int_{\Omega} |\phi|^p dz,
\end{equation}
where $\epsilon > 0$ and $C(\epsilon) > 0$. Substituting \eqref{L<} into \eqref{<L}, we obtain
\begin{equation}
\int_{\Omega} (\epsilon+(1-\epsilon)v(z))|\nabla_\gamma \phi|^p dz - (1-\epsilon) \int_{\Omega} \tilde{w}(z) |\phi|^p dz\geq
- (1-\epsilon) C(\epsilon) \int_{\Omega} |\phi|^p dz.
\end{equation}
The above inequality can be written as follows,
$$
\mathcal{R} := \frac{\int_{\Omega} (\epsilon+(1-\epsilon)v(z))|\nabla_\gamma \phi|^p dz - \int_{\Omega} (1-\epsilon) \tilde{w}(z) |\phi|^p dz}{\int_{\Omega} |\phi|^p dz} > - \tilde{C}(\epsilon),
$$
where $\tilde{C}(\epsilon) = (1-\epsilon) C(\epsilon)$. Hence,
$$
\inf_{0 \neq \phi \in C_0^{\infty}(\Omega \backslash\mathcal K)} \mathcal{R}> -\infty.
$$
This contradicts our assumption. The proof of Theorem \ref{parabolic eq} is now complete.
\end{proof}

\begin{proof}[Proof of Corollary \ref{cor1}]
To prove the Corollary \ref{cor1}, we need to show that 
\begin{equation}\label{rayleigh}
\inf_{0 \neq \phi \in C_c^{\infty}({B_{R}} \setminus (0,0))} \frac{A - B}{C} = -\infty,
\end{equation}
where $$A=\int_{B_{R}} \left(\epsilon+(1-\epsilon)\left(\log\frac{R}{\rho}\right)^{\alpha+p}\right)|\nabla_\gamma \phi|^p dz,$$ $$B=\int_{B_{R}} (1-\epsilon) c\left( \log \frac{R}{\rho} \right)^{\alpha} \frac{|x|^{\gamma p}}{\rho^{\gamma p+p}} |\phi|^p dz$$ and $$C=\int_{B_{R}} |\phi|^p dz.$$
By using \say{virtual} extremizer $\tilde{\phi}=\left(\log \frac{R}{\rho}\right)^{\frac{|\alpha+1|}{p}}$, we construct the test function $$\phi(\rho)=\psi_\delta(\rho)\left(\log \frac{R}{\rho}\right)^{\frac{|\alpha+1|}{p}+\epsilon}$$
for small $\delta>0$, $\epsilon>0$, and \( \psi_\delta \) is a smooth, radially symmetric function satisfying
\begin{equation}\label{psi}
    \psi_\delta(\rho) = \begin{cases}0, &\quad 0 \leq \rho \leq \delta, \\
    0 \leq \psi_\delta \leq 1, &\quad \delta \leq \rho \leq R-\delta, \\
    1, &\quad R-\delta\leq\rho \leq R,
    \end{cases}
\end{equation}
and $|\psi'_\delta(\rho)|\le \frac{C}{\delta}$.\\
A direct computation shows that 
\begin{equation}
|\nabla_\gamma\phi(\rho)|=
\begin{cases}    
0 & \text{if } 0<\rho<\delta,
\\
\frac{|x|^{\gamma}}{\rho^{\gamma}}\left|\frac{d}{d\rho}\left(\psi_\delta(\rho)\left(\log \frac{R}{\rho}\right)^{\frac{|\alpha+1|}{p}+\epsilon}\right)\right| & \text{if } \delta<\rho<R-\delta,\\
\left({\frac{|\alpha+1|}{p}+\epsilon}\right)\frac{|x|^{\gamma}}{\rho^{\gamma  +1}}\left(\log \frac{R}{\rho}\right)^{\frac{|\alpha+1|}{p}+\epsilon-1} & \text{if } R-\delta<\rho<R.
\end{cases}
\end{equation}\\
If $\phi:\Omega \to \mathbb{R}$ is radial, i.e., $\phi(z)=\phi(|z|)$, then
\begin{equation}\label{polar}
\int_{\Omega} \phi(|z|)\,dz
=
\Gamma
\int_{R_{1}}^{R_{2}} \rho^{Q-1}\phi(\rho)\,d\rho .
\end{equation}
Here
\[
\Gamma
=
\int_{a_{\theta}}^{b_{\theta}} d\theta
\int_{0}^{\pi} d\theta_{1}\cdots
\int_{0}^{\pi} d\theta_{k-2}
\int_{0}^{2\pi} d\theta_{k-1}
\int_{0}^{\pi} d\omega_{1}\cdots
\int_{0}^{\pi} d\omega_{d-2}
\int_{0}^{2\pi} d\omega_{d-1}\,
\Theta,
\]
and
\[
\Theta
=
\left(\frac{1}{1+\gamma}\right)^{k}
\lvert \sin\theta\rvert^{\frac{d}{1+\gamma}-1}
\cos^{k-1}\theta\,
\sin^{k-2}\theta_{1}\cdots
\sin\theta_{k-2}\,
\sin^{d-2}\omega_{1}\cdots
\sin\omega_{d-2}.
\]
Now, by using $|x|=\rho\,|\sin\theta|^{\frac{1}{1+\gamma}}$ and \eqref{polar}, we get
\begin{multline}
A=\int_{B_{R}} \left(\epsilon+(1-\epsilon)\left(\log\frac{R}{\rho}\right)^{\alpha+p}\right)|\nabla_\gamma \phi|^p dz\\= \mu\Bigg(\int_{\delta}^{R-\delta} \left(\epsilon+(1-\epsilon)\left(\log\frac{R}{\rho}\right)^{\alpha+p}\right)\left|\frac{d}{d\rho}\left(\psi_\delta(\rho)\left(\log \frac{R}{\rho}\right)^{\frac{|\alpha+1|}{p}+\epsilon}\right)\right|^p\rho^{Q-1}\,d\rho\\+\int_{R-\delta}^{R} \left(\epsilon+(1-\epsilon)\left(\log\frac{R}{\rho}\right)^{\alpha+p}\right)\left({\frac{|\alpha+1|}{p}+\epsilon}\right)^p\left(\log \frac{R}{\rho}\right)^{|\alpha+1|+\epsilon p-p}\rho^{Q-1-p}\,d\rho\Bigg)\\=\mu\Bigg(\int_{R-\delta}^{R} \left(\epsilon+(1-\epsilon)\left(\log\frac{R}{\rho}\right)^{\alpha+p}\right)\left({\frac{|\alpha+1|}{p}+\epsilon}\right)^p\left(\log \frac{R}{\rho}\right)^{|\alpha+1|+\epsilon p-p}\rho^{Q-1-p}\,d\rho\\+C_1\Bigg),
\end{multline}
where
\[
\mu
=
\Gamma
|\sin\theta|^{\frac{{\gamma} p}{1+\gamma}}.
\]
Next we can get
\begin{multline}
B=\int_{B_{R}} (1-\epsilon) c\left( \log \frac{R}{\rho} \right)^{\alpha} \frac{|x|^{\gamma p}}{\rho^{\gamma p+p}} |\phi|^p dz\\=\mu\Bigg(\int_\delta^{R-\delta} (1-\epsilon) c\left( \log \frac{R}{\rho} \right)^{\alpha} (\psi_\delta(\rho))^p\left(\log \frac{R}{\rho}\right)^{|\alpha+1|+\epsilon p} \rho^{Q-1-p} \,d\rho\\+\int_{R-\delta}^R (1-\epsilon) c\left( \log \frac{R}{\rho} \right)^{\alpha+|\alpha+1|+\epsilon p}    \rho^{Q-1-p} \,d\rho\Bigg)\\=\mu\Bigg(\int_{R-\delta}^R (1-\epsilon) c\left( \log \frac{R}{\rho} \right)^{\alpha+|\alpha+1|+\epsilon p}    \rho^{Q-1-p} \,d\rho+C_2\Bigg).
\end{multline}
In the same way we calculate
\begin{multline}
    C=\int_{B_{R}} |\phi|^p dz\\=\Gamma\left(\int_{\delta}^{R-\delta}(\psi_\delta(\rho))^p\left(\log \frac{R}{\rho}\right)^{|\alpha+1|+\epsilon p}\rho^{Q-1}\,d\rho+\int_{R-\delta}^R\left(\log \frac{R}{\rho}\right)^{|\alpha+1|+\epsilon p}\rho^{Q-1}\,d\rho\right),
\end{multline}
which is equal to some constant, if $\epsilon\xrightarrow{}0^+$.
Lastly, if we calculate $A-B$, as $\epsilon\xrightarrow{}0^+$, we get
\begin{multline}
    A-B= \\\mu\Bigg(\int_{R-\delta}^{R} \left(\epsilon+(1-\epsilon)\left(\log\frac{R}{\rho}\right)^{\alpha+p}\right)\left({\frac{|\alpha+1|}{p}+\epsilon}\right)^p\left(\log \frac{R}{\rho}\right)^{|\alpha+1|+\epsilon p-p}\rho^{Q-1-p}\,d\rho\\-\int_{R-\delta}^R (1-\epsilon) c\left( \log \frac{R}{\rho} \right)^{\alpha+|\alpha+1|+\epsilon p}    \rho^{Q-1-p} \,d\rho+C_1-C_2\Bigg)\\=\mu\Bigg(\left({\frac{|\alpha+1|^p}{p^p}-c}\right)\int_{R-\delta}^{R} \left(\log\frac{R}{\rho}\right)^{-1}\rho^{Q-1-p}\,d\rho+C_3\Bigg),
\end{multline}
where $$\frac{|\alpha+1|^p}{p^p}-c<0,$$
$$\int_{R-\delta}^{R} \left(\log\frac{R}{\rho}\right)^{-1}\rho^{Q-1-p}\,d\rho=+\infty.$$
Thus, the infimum of the Rayleigh quotient is
\begin{equation}
\inf_{0 \neq \phi \in C_c^{\infty}({B_{R}} \setminus (0,0))} \frac{A - B}{C} = -\infty,
\end{equation}
which proves the Corollary.
\end{proof}

\begin{proof}[Proof of Corollary \ref{cor2}]
To prove Corollary \ref{cor2}, we need to show that
$$
\inf_{0 \neq \phi \in C_c^{\infty}(B_R)} \frac{\int_{B_R} |\nabla_\gamma \phi|^p dz - c(1-\epsilon) \int_{B_R} \frac{|x|^{\gamma p}}{(R-\rho)^{p}\rho^{\gamma p}} |\phi|^p dz}{\int_{B_R} |\phi|^p dz} = -\infty.
$$
By using \say{virtual} extremizer $\tilde{\phi}=(R-\rho)^{\frac{p-1}{p}}$, we construct the test function 
$$\phi(\rho)=\psi_\delta(\rho)(R-\rho)^{\frac{p-1}{p}+\epsilon},$$
where $\psi_\delta$ is a smooth, radially symmetric function as defined in \eqref{psi}.
A direct computation shows that 
\begin{equation}
|\nabla_\gamma\phi(\rho)|=
\begin{cases}    
0 & \text{if } 0<\rho<\delta,
\\
\frac{|x|^\gamma}{\rho^\gamma}\left|\frac{d}{d\rho}\left(\psi_\delta(\rho)(R-\rho)^{\frac{p-1}{p}+\epsilon}\right)\right| & \text{if } \delta<\rho<R-\delta,\\
\left({\frac{p-1}{p}+\epsilon}\right)(R-\rho)^{\frac{p-1}{p}+\epsilon-1} & \text{if } R-\delta<\rho<R.
\end{cases}
\end{equation}
Now, by using $|x|=\rho\,|\sin\theta|^{\frac{1}{1+\gamma}}$ and \eqref{polar}, we get

\begin{multline}
\int_{B_{R}} |\nabla_\gamma \phi|^p dz= \mu\Bigg(\int_{\delta}^{R-\delta} \left|\frac{d}{d\rho}\left(\psi_\delta(\rho)(R-\rho)^{\frac{p-1}{p}+\epsilon}\right)\right|^p\rho^{Q-1}\,d\rho\\+\int_{R-\delta}^{R} \left({\frac{p-1}{p}+\epsilon}\right)^p(R-\rho)^{p-1+\epsilon p -p}\rho^{Q-1}\,d\rho\Bigg)\\=\mu\Bigg(\int_{R-\delta}^{R} \left({\frac{p-1}{p}+\epsilon}\right)^p(R-\rho)^{-1+\epsilon p}\rho^{Q-1}\,d\rho+C_1\Bigg),
\end{multline}
Next we get
\begin{multline}
c(1-\epsilon) \int_{B_R} \frac{|x|^{\gamma p}}{(R-\rho)^{p}\rho^{\gamma p}} |\phi|^p dz\\=\mu\Bigg(c(1-\epsilon)\int_\delta^{R-\delta}(R-\rho)^{-p}(\psi_\delta(\rho))^p(R-\rho)^{p-1+\epsilon p} \rho^{Q-1} \,d\rho\\+c(1-\epsilon)\int_{R-\delta}^R(R-\rho)^{-p}(R-\rho)^{p-1+\epsilon p} \rho^{Q-1} \,d\rho\Bigg)\\=\mu\Bigg(c(1-\epsilon)\int_{R-\delta}^R (R-\rho)^{-1+\epsilon p} \rho^{Q-1} \,d\rho + C_2\Bigg).
\end{multline}
In the same way, we calculate
\begin{multline}
    \int_{B_{R}} |\phi|^p dz\\=\Gamma\left(\int_{\delta}^{R-\delta}(\psi_\delta(\rho))^p(R-\rho)^{p-1+\epsilon p}\rho^{Q-1}\,d\rho+\int_{R-\delta}^R(R-\rho)^{p-1+\epsilon p}\rho^{Q-1}\,d\rho\right),
\end{multline}
which is equal to some constant, if $\epsilon\xrightarrow{}0^+$.

Lastly, we evaluate the infimum of the Rayleigh quotient by taking the limit as $\epsilon \to 0^+$. Combining the previous estimates, we have:

\begin{multline}
    \inf_{0 \neq \phi \in C_c^{\infty}(B_R)} \frac{\int_{B_R} |\nabla_\gamma \phi|^p dz - c(1-\epsilon) \int_{B_R} \frac{|x|^{\gamma p}}{(R-\rho)^{p}\rho^{\gamma p}} |\phi|^p dz}{\int_{B_R} |\phi|^p dz} \\=\inf_{0 \neq \phi \in C_c^{\infty}(B_R)}\frac{\mu\Bigg(\left(\left({\frac{p-1}{p}+\epsilon}\right)^p-c(1-\epsilon)\right)\int_{R-\delta}^{R} (R-\rho)^{-1+\epsilon p}\rho^{Q-1}\,d\rho+C_3\Bigg)}{C_4}\\=\frac{\mu\Bigg(\left(\left({\frac{p-1}{p}}\right)^p-c\right)\int_{R-\delta}^{R} (R-\rho)^{-1}\rho^{Q-1}\,d\rho+C_3\Bigg)}{C_4}=-\infty, 
\end{multline}
as $$\left(\left({\frac{p-1}{p}}\right)^p-c\right)<0$$
and 
$$\int_{R-\delta}^{R} (R-\rho)^{-1}\rho^{Q-1}\,d\rho=+\infty.$$ 
This proves the Corollary.
\end{proof}

\bibliographystyle{alpha}
\bibliography{citation}

\end{document}